\documentclass[3p,times]{elsarticle}

\usepackage{amssymb}
\usepackage[figuresright]{rotating}

\newdefinition{remark}{Remark}
\newproof{proof}{Proof}

\usepackage{graphicx}
\usepackage{amsmath}
\usepackage{exscale,times}
\usepackage{float}
\usepackage{color}
\usepackage{setspace}
\usepackage{amssymb}
\usepackage{graphicx}
\usepackage{amsmath}
\usepackage{exscale,times}
\usepackage{epsfig}
\usepackage{subfig}
\usepackage{float}
\usepackage{color}
\usepackage{setspace}
\usepackage{mathtools}
\allowdisplaybreaks               %% denklem 2. sayfaya tasarsa, ayarl�yor.
\begin{document}
\doublespacing
\begin{frontmatter}



%% e.g. \dochead{17th International Conference on Dynamical Processes in Excited States of Solids}

\title{APPLICATION OF ADAPTED-BUBBLES TO THE HELMHOLTZ EQUATION WITH LARGE WAVE NUMBERS IN 2D}


%\author{}
%
%\address{}

%\author[rvt]{Adem Kaya}
%\ead{kaya@uni-potsdam.de}

%\author[rvt]{Melina Freitag} 
%\ead{melina.freitag@uni-potsdam.de}

\author[rvt]{Adem Kaya\corref{cor}}  
\cortext[cor]{Corresponding author} 
\ead{kaya@uni-potsdam.de}

%\author[rvt]{}% <- Another address for this author
%\ead{} 

\address[rvt]{Institut f\"{u}r Mathematik, Universit\"{a}t Potsdam
Karl-Liebknecht-Str. 24-25
14476 Potsdam/Golm
Germany}

\begin{abstract}
  An adapted bubble approach  which is a modifiation of the residual-free bubbles (RFB) method, is proposed for the Helmhotz problem in 2D.  A  new two-level finite element method is introduced for the approximations of the bubble functions.  Unlike the other equations such as the advection-diffusion equation, RFB method when applied to the Helmholtz equation, does not depend on another stabilized method  to obtain  approximations to the solutions of the sub-problems.  Adapted bubbles (AB) are obtained by a simple modification of the sub-problems. This modification increases the accuracy of the numerical solution  impressively. The AB method is able to solve the Helmholtz equation efficiently in 2D up to $ch=3.5$ where $c$ is the wave number and $h$ is the mesh size.  We provide  analysis to show how the AB method  mitigates the pollution error.          
\end{abstract}

\begin{keyword}
  
Helmholtz equation, adapted bubbles, residual-free bubbles,  two-level finite element  

\end{keyword}

\end{frontmatter}

%%\linenumbers


\section{ Introduction}

Enriching linear finite element space with \textit{residual-free bubble functions} is a general framework for the discretizations of the problems\cite{RFBhelmholtz, choosingbuubles,nesliturkcdr, nesliturktez,unlocking, masslumping}.     
These functions strongly satisfy the original differential equations and hence obtaining the bubble functions is generally as difficult as solving the original problem such as  the convection-diffusion equation \cite{choosingbuubles}.  Unlike it was stated in \cite{RFBhelmholtztwolevel}, we will show that this is not the case for the Helmholtz problem. Obtaining the bubble functions on triangular elements is easier  than solving the original problem. The standard Galerkin finite element method can be used with a coarse mesh to obtain efficient approximations to the bubble functions.

The residual-free bubbles method produces the exact solution of linear dierential equations in the one-dimensional case. However, the method in higher-dimensions  is approximate and as we will show for the Helmholtz problem in this article, its contribution to the stabilization of the standard Galerkin method is very poor. We modify the residual-free bubbles (RFB) method in 2D by multiplying the right-hand side of the bubble equations with a constant. This operation impressively increase the accuracy. The new bubbles are no more residual-free and we call them \textit{ adapted bubbles} (AB). We provide the optimal values of the constants for the triangular and rectangular elements separately. We apply a two-level finite element method using the standard  Galerkin finite element method to get approximations to the bubble functions. 

We provide  analysis to show how  the AB method  mitigates the pollution error. To this end. we approximate the bubble functions with piece-wise defined linear functions so-called \textit{pseudo}-bubbles. The analysis give rise to a fourth order finite difference scheme with seven-point stencil for plane waves. It is perfectly applicable in polygonal and triangular domains. We use this method to do comparison with the AB method.

Standard discretizations when applied to the Helmholtz problem  suffer from the pollution effect when the wave number is large \cite{Babuska1997}. Moreover standard iterative solvers are ineffective in obtaining the solutions of the discrete Helmholtz equation \cite{Ernst2012}. There is a great effort in literature to overcome these difficulties. Among the discretization techniques, there are 
finite difference \cite{Singer1998,Feng2011}, finite element \cite{Ihlenburg1995,Strouboulis2006,Alvarez2006}, discontinuous
Galerkin \cite{paul, Haijun}, virtual element  \cite{Perugia2016}, and boundary element methods \cite{Kirkup2019}. At the same time, there is a great effort to develop efficient preconditioners, such as multigrid \cite{Elman2002, Erlangga2006} and domain decomposition methods \cite{Gander2003, Stolk2013}. The AB method proposed in this article, does not suffer from the pollution effect for very large wave numbers. It is by far superior than the fourth order accurate scheme proposed here.

The rest of this paper is organized as follows. We review the RFB method  for the Helmholtz equation  in Section \ref{sec:section1}. We explain how to implement two-level finite element method in 1D and provide analysis to show the contribution of the  bubble functions in reducing  the pollution error in Section \ref{sec:section2}. Section \ref{sec:section3} is devoted for the  analysis of the RFB method in 2D.
We propose the AB method for triangular elements in Section  \ref{sec:section4}. Several numerical experiments are provided in Section \ref{sec:section5}. The AB method is considered with rectangular elements in Section \ref{sec:section6}.  
We finish with concluding remarks in Section \ref{sec:section7}.

\section{Residual-free bubbles method (RFB) for the Helmholtz equation}\label{sec:section1}

We start with considering the Helmholtz problem in 1D with Dirichlet boundary conditions on unit interval.
\begin{eqnarray}
 \Bigg\{
 \begin{array}{ll}
 -u^{\prime \prime} -c^2u(x) = f(x), \qquad x\in I\\ 
 u(0) = 0, \qquad u(1) = 1,
 \end{array}
 \label{eqn:Helmholtz1dd}
\end{eqnarray}
where we assume that the wave number $c$ is constant.  
RFB method can be summarized as follows. Let's start 
with recalling abstract variational formulation of 
(\ref{eqn:Helmholtz1dd}): Find $u \in H^1(I)$ such that 
\begin{equation*}
  a(u,v)=(f,v), \qquad \forall v \in H^1(I),
\end{equation*}
where
\begin{equation*}
  a(u,v) =  \int_I u' v' dx - c^2\int_I  u v dx
\end{equation*}
and 
\begin{equation*}
 (f,v) = \int_I  f v dx.
\end{equation*}
Define $V_h \subset H^1(I)$ as a finite-dimensional space. Then the Galerkin finite element method reads: Find $u_h \in V_h$ such that
\begin{equation*}
  a(u_h,v_h)=(f,v_h), \qquad \forall v_h \in V_h.
\end{equation*}
We now decompose the space $V_h$ as $V_h=V_L \bigoplus V_B$, where $V_L$ is the space of continuous piecewise linear polynomials and $V_B=\bigoplus_K B_K$ with
$B_K=H_0^1(K)$. From this decomposition, every $v_h \in V_h$ can be written in the form $v_h=v_L+v_B$, where $v_L \in V_L$ and $v_B \in V_B$. 
Bubble component $u_B$ of $u_h$ satisfy the original differential equation in an element $K$ strongly, i.e.
\begin{equation}
 \mathcal{L}u_B = -\mathcal{L}u_L + f \qquad \text{in} \qquad K,
 \label{eqn:bubble1}
\end{equation}
subject to boundary condition,
\begin{equation}
 u_B=0 \qquad \text{on} \qquad \partial K. 
  \label{eqn:bubble2}
\end{equation}
Since the support of bubble $u_B$ is contained within the element $K$, we can make a static condensation for the bubble part, getting directly the $V_L$-
projection $u_L$ of the solution $u_h$ \cite{RFBhelmholtz}. This can be done as follows. Using $V_h=V_L \bigoplus V_B$, the finite element approximation 
reads: Find $u_h=u_L+u_B$
 in $V_h$ such that 
\begin{equation}
 a(u_L,v_L)+a(u_B,v_L)=(f,v_L), \qquad \forall v_L \in V_L.
 \label{eqn:varfordecom}
\end{equation}

\section{A two level finite element method} \label{sec:section2}
  In order to find bubble part $u_B$ of the solution,  we need to solve (\ref{eqn:bubble1})-(\ref{eqn:bubble2}). The problems defined by equations (\ref{eqn:bubble1})-(\ref{eqn:bubble2}) is addressed by solving instead
 \begin{eqnarray}
\left\lbrace
\begin{array}{ll}
  - \varphi_i''  -c^2 \varphi_i =  c^2 \psi_i \quad \text{in} \quad K, \quad (i=1,...,n_{en})\\
 \varphi_i=0 \quad \text{on} \quad \partial K,
 \end{array}\right.
\label{eqn:bubble3} 
\end{eqnarray} 
and
 \begin{eqnarray}
\left\lbrace
\begin{array}{ll}
 -  \varphi_f'' -c^2  \varphi_f = f \quad \text{in} \quad K,\\
  \varphi_f = 0 \quad \text{on} \quad \partial K.
 \end{array}\right.
\label{eqn:bubble4} 
\end{eqnarray}
where $n_{en}$ is the number of element nodes.
Thus if 
\begin{equation}
 u_L=\sum\limits_{i} d_i^K \psi_i,
 \label{eqn:ul}
\end{equation}
then
\begin{equation}
 u_B=\sum\limits_{i} d_i^K \varphi_i +\varphi_f.
 \label{eqn:ub}
\end{equation}
Substituting (\ref{eqn:ul}) and (\ref{eqn:ub}) into (\ref{eqn:varfordecom}), we get the matrix formulation 
\begin{equation}
\sum_{K}\sum\limits_{i}^{n_{en}} d_i^K\left( (\psi^{'}_i,\psi_j^{'}) -c^2  
(\psi^{}_i,\psi_j^{})  -c^2(\varphi_i^{},\psi_j^{})  )=(f,\psi_j) \right)
\end{equation}
at the global level where $d_i$ are the finite element approximations to the solutions at the nodes.

Numerical solution of the bubble problems (\ref{eqn:bubble3}) and (\ref{eqn:bubble4})  generally requires using a nonstandard method  such as for the case of the advection-diffusion equation. This  makes the RFB method  dependent on another stabilized method when applying it as a two-level finite element method. In \cite{RFBhelmholtztwolevel},   the Galerkin-least-squares method (GLS) was used to get approximations to the bubble functions  in solving the  the Helmholtz equation.  Although this is true for the advection-diffusion equation, there is no need to use a nonstandard method to get approximations to the bubble functions when the Helmholtz problem is under consideration. We explain this this fact in 1D. It is well known that standard discretizations are effective up to $ch=0.6$. Suppose that we have a discretization of the domain such that $ch=0.6$. Even if we use 3 nodes on the sub-domain (element), $ch_e$ ($h_e$ is mesh size on the sub-domain) becomes less than $0.6$. If $ch=3$ on the global mesh, then using 11 nodes for the sub-problems makes $ch_e=0.3$. More precisely, it is always true that $ch_e <ch$.

It is true that GLS computation is known to incur at most marginal increase in computational cost over the standard Galerkin method. However, GLS for the sub-problems  may lead to misinterpretations related to the bubble functions.

Another way of obtaining the bubble function is to use separation of variables when rectangular elements are used \cite{RFBhelmholtz}. However, this gives rise to a series solution of the bubble function for which it must be truncated. For a good accuracy,  200 terms are used in  \cite{RFBhelmholtz} which is  computationally not so effective. Another drawback is that this approach is limited to the rectangular elements.

\subsection{Analysis of the pollution effect of the sub-problems }
The inequality $ch_e <ch$ is an indication that the sub-problems are easier to solve; however, we must analyse the pollution effect for the sub-problems for large wave numbers. It is well known that the condition 
 $c^2h<1$ is sufficient to guarantee that the error of the Galerkin solution is of the same magnitude as the error of the best approximation \cite{Babuska1997}. This condition is necessary when the size the domain is fixed for increasing $c$.  More precisely, the exact solution is very oscillatory. In our case, the exact solutions of the sub-problems are not oscillatory. When $ch<\pi$, the exact solutions of the sub-problems are always in the form of a half wave as the homogenous Dirichlet boundary condtion is applied everywhere on the boundary. When $\pi<ch<2\pi$, the exact solutions of the sub-problems are always in the form of  a single wave.  In this regime, the standard Galerkin method is pollution free for the sub-problems for any wave number.   
 Note that, in simulations, 12 nodes per wave is generally chosen which correponds to $ch\approx0.62$.

 We use the  standard Galerkin finite element method with piecewise linear basis functions to approximate the bubble functions. Note that  the bubble problems (\ref{eqn:bubble3}) and (\ref{eqn:bubble4}) can be solved independently and hence  parallel processors can be used to carry out these computations efficiently. When uniform meshes are used and the right hand side function $f$ is constant, construction of the system matrix is as cheap as construction of the system matrix of the standard Galerkin finite element method.

 \subsection{Shape of the bubble functions and pseudo-bubbles}

We have shown that the RFB method is not dependent on another stabilized method to get approximations to the bubble functions when a two-level finite element method is used due to the non-oscillatory behavior of the exact solutions of the sub-problems. This non-oscillatory behavior of the exact solutions opens a gateway to approximate these bubble functions with piecewise-defined linear simple functions. These approximations are called \textit{pseudo-bubbles} and constructed considering the shape of the bubble functions. Pseudo-bubbles were applied to the advection-diffusion-reaction equation in \cite{Asendur,Asendur2}.        

\begin{figure}[H]
\center
 \includegraphics[width=7.5cm]{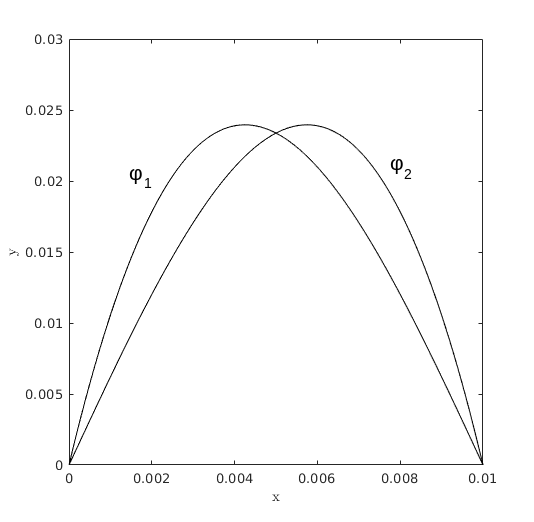}
 \includegraphics[width=8.7cm]{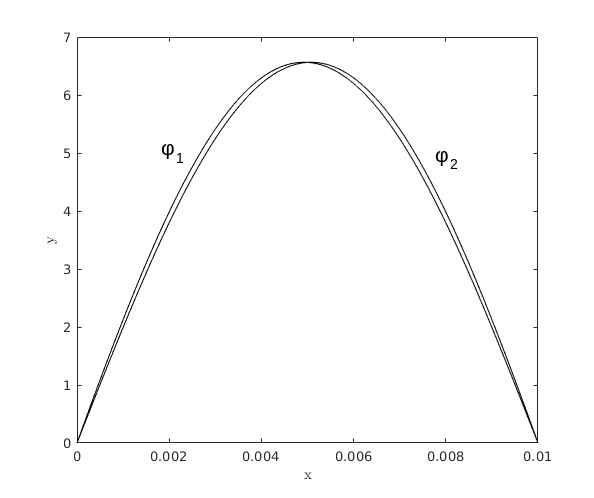}   
\caption{Bubble functions for $c=60$ (left) and $c=300$ (right).}
\label{fig:1dbubbles}
\end{figure}

\begin{figure}[H]
 \includegraphics[width=8cm]{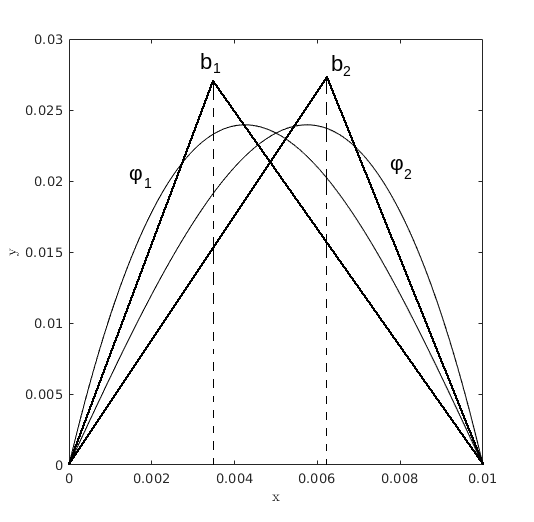}
  \includegraphics[width=8cm]{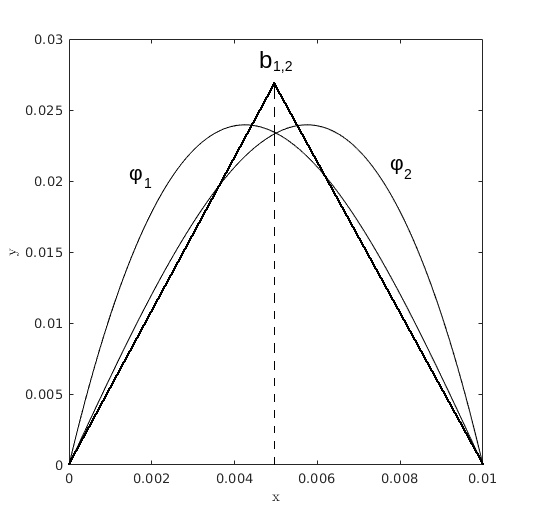}
\caption{Two different approaches to approximate the bubble functions.}
\label{fig:1dpseudobubbles}
\end{figure}

Here, we consider the case $ch<\pi$ for which the bubble functions are in the form of a half wave. We present the bubble functions $\varphi_{1,2}$ in Figure \ref{fig:1dbubbles} for $c=60,300$ when  $h=0.01$. Efficient yet cheap approximations to these bubble functions with piecewise-defined linear functions are given in Figure \ref{fig:1dpseudobubbles}. While on the left in Figure \ref{fig:1dpseudobubbles}, two pseudo-bubbles are used, it is possible to approximate $\varphi_{1,2}$ with a single pseudo-bubble. The humps of the bubble functions  $\varphi_{1,2}$ come closer as $ch$ increases (see Figure \ref{fig:1dbubbles}. Thus, for larger $ch$, we can derive more efficient approximations to the bubble functions when  a single pseudo-bubble is used. The advantage of using a single pseudo-bubble is that the maximum of the pseudo-bubble occurs in the middle of the element. Applying the minimization technique applied in \cite{Asendur}, one can find the optimal heights  and locations of the peaks of the pseudo-bubbles. When a single pseudo-bubble is used, one can easily calculate integrals in the finite element formulation. This will be important in modifying the RFB method in 2D.

\subsection{Analysis of the  pseudo-bubbles}
In order to see how the residual-free bubbles method  overcome the pollution effect, we first consider linear finite element method for (\ref{eqn:Helmholtz1dd}) when $f(x)=0$, $u(0)=0$ and $u(1)=\sin(c)$.  The exact solution is $u(x)=\sin(cx)$. We obtain the truncation error and see how the error deteriorates as $c$ increases for fixed $ch$ which is simply the pollution effect. To this end, we jump to the finite difference equivalance of the linear finite method. Let $U_j$ represents the numerical solution and choose $n$  equally distributed nodes for which $h=1/(n-1)$. Taking the integrals in linear finite element formulation and scaling by $h$  gives
\begin{equation}
-\frac{U_{j+1} -2U_{j} + U_{j-1}}{h^2} - c^2 \frac{U_{j+1} +4 U_{j} + U_{j-1}}{6} =0, \qquad j=2,...,n-2.  
\label{eqn:finitdifference1d}
\end{equation}
From the boundary conditions,
\begin{equation*}
 U_{1} = 0, \qquad U_{n} = \sin(c).
\end{equation*}
Using the Taylor series expansion, we get the truncation error for (\ref{eqn:finitdifference1d}).
\begin{equation}
\tau(x) = -\frac{c^2 h^2}{6} u^{\prime \prime} - \left(\frac{c^2h^4}{72} + \frac{h^2}{12} \right)u^{(4)} - \frac{h^4}{360} u^{(6)} + \mathcal{O}(h^6).    
\label{eqn:truncationerror1d}
\end{equation}
The pollution effect can not be seen from (\ref{eqn:truncationerror1d}).
To see it,  we substitute the exact solution   $u(x)=\sin(cx)$ into (\ref{eqn:truncationerror1d}).  
\begin{equation*}
\tau(x) = \frac{c^4 h^2}{6} \sin(cx) - \left(\frac{c^6h^4}{72} + \frac{c^4h^2}{12} \right) \sin(cx) + \frac{c^6 h^4}{360} \sin(cx) + \mathcal{O}(c^8h^6).    
\label{eqn:truncationerror1d2}
\end{equation*}
Rearranging   the above equation gives
\begin{equation}
\tau(x) = \sin(cx) \left( \frac{c^4 h^2}{12}  - \frac{c^6h^4}{90}  + \mathcal{O}(c^8h^6)\right).    
\label{eqn:truncationerror1d2}
\end{equation}

When the exact solution is oscillatory, that is, $c$ is large, the  term  $c^4h^2/12$ in (\ref{eqn:truncationerror1d2})  becomes large, that is, $\tau(x)$ is large, even if $ch=\text{constant}$ is small. This is called the pollution effect. When $ch$ is sufficiently small, there is no phase error for the Dirichlet problem when $c$ is large, however  when a Neumann or Robin boundary condition is used, phase error is also observed.

The simplest way to mitigate this pollustion effect is to choose $c^2h$ sufficiently small. However, this requires intractable matrices in higher dimensions.      
The general idea in literature is to decrease the effect of  the first  few terms in (\ref{eqn:truncationerror1d2}) so that tractable matrix sizes can be obtained. For example, using higher order accurate methods of finite difference  or higher order polynomials  finite element may allow to eliminate the first few terms. If the first term can be eliminated, then the requirement to mitigate the pollution effect  reduces to $c^{3/2}h$ being sufficiently small. However, higher order methods generally use more points and this  increases the nonzero entries of the  matrices.        

In order to get a deeper insight of the working principle of the residual-free bubbles method in mitigating the pollution effect, we consider the pseudo-bubbles in Figure \ref{fig:1dpseudobubbles} on the right. This choice allows us to take the integrals containing the bubble functions,  explicitly.   
\begin{figure}[H]
\center 
 \includegraphics[width=8cm]{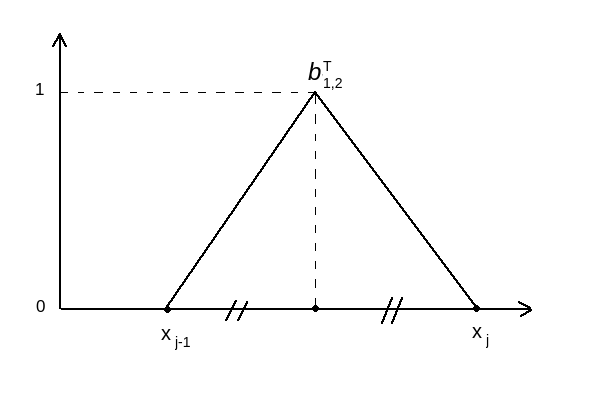}
\caption{Basis functions employed in the approximation of bubble functions.}
\label{fig:bubblesbasis} 
\end{figure} 
We can define the pseudo-bubbles  $b_{1,2}$ using the basis functions $b_{1,2}^T$ represented in Figure \ref{fig:bubblesbasis} and heights of $b_{1,2}$, i.e.,  $\alpha_{1,2}$. More precisely,
\begin{equation}
 b_{1,2} = \alpha_{1,2}b_{1,2}^T,
\label{eqn:bubbledefinition}
\end{equation}
where $\alpha_1=\alpha_2$.  Applying the  technique proposed in \cite{Asendur} (set $\xi=h/2, \epsilon=1, \sigma=-c^2$ in equation (13) in \cite{Asendur}) gives
\begin{equation}
 \alpha_1 = \frac{3c^2h^2}{4(12-c^2h^2)}.
 \label{eqn:bubblealf}
\end{equation}

Taking the integrals in  (\ref{eqn:varfordecom}) making use of (\ref{eqn:bubbledefinition})  gives the finite difference formula 
\begin{equation}
-\frac{U_{j+1} -2U_{j} + U_{j-1}}{h^2} - c^2 \frac{U_{j+1} +4 U_{j} + U_{j-1}}{6} - \alpha_1 c^2 \frac{U_{j+1} +2 U_{j} + U_{j-1}}{4} =0, \qquad j=2,...,n-2.  
\label{eqn:pseudobubblefinite}
\end{equation}
Using   Taylor  expansions of $u(x\pm h)$, definition of $\alpha_1$ given in (\ref{eqn:bubblealf}) and the exact solution $u(x) = \sin(cx)$, we obtain the truncation error
\begin{equation*}
\tau(x) = \sin(cx) \left( \frac{c^4 h^2}{12}  - \frac{c^6h^4}{90}  + \mathcal{O}(c^8h^6)\right) + \sin(cx) \left( -\frac{3c^4 h^2}{4(12-c^2h^2)}  + \frac{3c^6h^4}{16(12-c^2h^2)}  + \mathcal{O}(c^8h^6)\right) .    
\end{equation*}
Rearranging the right-hand side of the above formula we end up with
\begin{equation}
\tau(x) = \sin(cx) \left( c^4h^2 \left( \frac{1}{12} -\frac{3}{4(12-c^2h^2)} \right)  + c^6h^4 \left( - \frac{1}{90}  + \frac{3}{16(12-c^2h^2)} \right) + \mathcal{O}(c^8h^6)\right) .    
\label{eqn:truncationerror1d5}
\end{equation}
While the coefficient of $c^4h^2$ in (\ref{eqn:truncationerror1d2}) is $1/12$, the coefficient of $c^6h^4$ is $1/90$ in magnitude. When $ch<\sqrt{15/2} \approx 2.73$,   
the coefficient of $c^4h^2$ in (\ref{eqn:truncationerror1d5})  
is smaller in magnitude and it is close to zero when $ch\approx \sqrt{3}$.   
Moreover, when $ch<\sqrt{57/16}\approx 1.88$, the coefficient of $c^6h^4$ in (\ref{eqn:truncationerror1d5}) becomes smaller in magnitude. 
The approximate bubbles shows how  the pollution effect is reduced. It is known that RFB method for 1D linear equations is exact \cite{RFBhelmholtztwolevel}. This means that it automatically makes the coefficient of all powers $c^{n+2}h^n$, $n=2,3,...$, zero. A good approximation to the residual-free bubbles significantly reduce the pollution effect. The RFB method is approximate in 2D. It is well known that the contributions of the residual-free bubble functions to the stabilization of the Galerkin method is very poor.   The observations we made here will be helpful to further increase the accuracy of the method in 2D. We will modify the sub-problems in 2D and use adapted bubbles to further increase the accuracy of the bubble approach.        

\section{The RFB method   in 2D}
\label{sec:section3}

\begin{figure}[H]
\center
 \includegraphics[width=8cm]{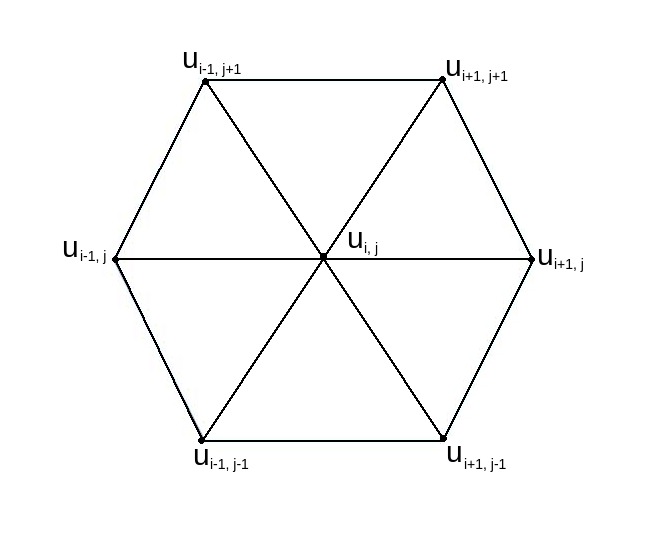} 
\caption{Equilateral triangles surrounding a node.}
\label{fig:2dstencil}
\end{figure}

We have shown that the RFB method is able to solve the Helmholtz problem in 1D cheaply and efficiently for very large wave numbers. As it was stated in \cite{RFBhelmholtztwolevel}, RFB method is not as  efficient in 2D as in 1D. To show this fact, we consider the following problem on an equilateral triangular shaped domain with vertices $(0,0)$, $(0,1)$ and $(0.5, \sqrt{3}/2)$. 
\begin{eqnarray}
\left\lbrace
 \begin{array}{ll}
 -\Delta u -c^2u = 0, \quad \text{in} \quad \Omega, \\
  u(x,y) = \sin(cx\sin(\theta) + cy\cos(\theta)),\quad \text{on} \quad \partial \Omega_D,
\end{array} \right.
 \label{eqn:Helmholtz2ddirichlet}
\end{eqnarray}
where the exact solution is $u(x,y)=\sin(cx\sin(\theta) + cy\cos(\theta))$. We use equilateral triangular elements with linear basis functions.    
We solve the following equations  on element level to get the RFB functions.
\begin{eqnarray}
\left\lbrace
\begin{array}{ll}
  - \Delta \varphi_i  -c^2 \varphi_i =   c^2 \psi_i \quad \text{in} \quad K, \quad (i=1,...,n_{en})\\
 \varphi_i=0 \quad \text{on} \quad \partial K,
 \end{array}\right.
\label{eqn:2dbubble3} 
\end{eqnarray} 
 and
 \begin{eqnarray}
\left\lbrace
\begin{array}{ll}
 - \Delta \varphi_f -c^2  \varphi_f = f \quad \text{in} \quad K,\\
  \varphi_f = 0 \quad \text{on} \quad \partial K,
 \end{array}\right.
\label{eqn:2dbubble4} 
\end{eqnarray} 
where $n_{en}=3$. $\varphi_i$ and $ \psi_i$, $(i=1,2,3)$ are the RFB  and the linear basis functions, respectively.
Linear   finite element method with a coarse mesh can be used to obtain efficient approximations to the bubble functions. To do some analyses, we  approximate   the RFB functions with piecewise-defined linear functions with the maximum  at the centroid of the element. Let $b_{1,2,3}=\alpha_2b^T$ be the approximation to the bubble functions where $b^T$ is the linear basis bubble function that assumes zero at the vertices of the element and one at the centroid of the element. Applying the same procedure we applied in 1D (see \cite{Asendur2} for more details), gives
\begin{equation}
 \alpha_2^2 \int_K \nabla b^T \nabla b^T dS - \alpha_2^2 c^2 \int_K b^T b^T dS = \alpha_2 \int_K \psi_i b^T dS, \qquad i=1,2,3.
\end{equation}
Solving the above equation for $\alpha_2$ and calculating the integrals for $\psi_1$ gives 
\begin{equation}
 \alpha_2 = \frac{2c^2h^2}{3(72-c^2h^2)}.
 \label{eqn:alfagalerkin}
\end{equation}
Considering 6 adjacent elements as shown in Figure \ref{fig:2dstencil}, the RFB method is equivalent to the following finite difference scheme.  

\begin{eqnarray}
\begin{array}{ll}
   \frac{6 U_{i,j} - U_{i-1,j} - U_{i-1,j+1} - U_{i+1,j+1} - U_{i+1,j} - U_{i+1,j-1} -U_{i-1,j-1}}{\sqrt{3}h^2} 
    -c^2  \frac{6 U_{i,j} + U_{i-1,j} + U_{i-1,j+1} + U_{i+1,j+1} + U_{i+1,j} + U_{i+1,j-1} + U_{i-1,j-1}}{8\sqrt{3}} \\
 -\alpha_2 c^2  \frac{3 U_{i,j} + U_{i-1,j} + U_{i-1,j+1} + U_{i+1,j+1} + U_{i+1,j} + U_{i+1,j-1} + U_{i-1,j-1}}{6\sqrt{3}}, \qquad i,j=1,...,N_{\text{int}},  
\end{array}
 \label{eqn:2dfinitedifference}
\end{eqnarray}
where
$U_{i,j} \approx u(x,y)$,  $U_{i-1,j} \approx u(x-h,y)$, $U_{i-1,j+1} \approx u(x-\frac{h}{2},y + \frac{\sqrt{3}}{2}h)$, 
$U_{i+1,j+1} \approx u(x+\frac{h}{2},y + \frac{\sqrt{3}}{2}h)$,
$U_{i+1,j} \approx u(x+h,y)$,
$U_{i+1,j-1} \approx u(x+\frac{h}{2},y - \frac{\sqrt{3}}{2}h)$ and
$U_{i-1,j-1} \approx u(x-\frac{h}{2},y - \frac{\sqrt{3}}{2}h)$.
Note that when $\alpha_2=0$, (\ref{eqn:2dfinitedifference}) is equivalent to the linear finite element method with equilateral triangular element. 

To analyse the RFB method, we substitute  the Taylor expansions of the exact solution at the grid points. The derivation of the truncation error is given in (\ref{eqn:truncationerror2d})-(\ref{eqn:truncationerror2d3}). In our analysis, we will examine the coefficients of $c^4h^2$ and $c^6h^4$ in (\ref{eqn:truncationerror2d3}), that is, $C_2$ and $C_1$ in (\ref{eqn:twocoefficient}).

\begin{sidewaysfigure}

\begin{eqnarray}
 \begin{array}{ll}
\tau(x,y) =  \frac{6 u(x,y) - u(x-h,y) - u(x-\frac{h}{2},y + \frac{\sqrt{3}}{2}h) - u(x+\frac{h}{2},y + \frac{\sqrt{3}}{2}h) - u(x+h,y) - u(x+\frac{h}{2},y - \frac{\sqrt{3}}{2}h) -u(x-\frac{h}{2},y - \frac{\sqrt{3}}{2}h)}{\sqrt{3}h^2}  -c^2  \frac{6 u(x,y) + u(x-h,y) + u(x-\frac{h}{2},y + \frac{\sqrt{3}}{2}h) + u(x+\frac{h}{2},y + \frac{\sqrt{3}}{2}h) + u(x+h,y) + u(x+\frac{h}{2},y - \frac{\sqrt{3}}{2}h) +u(x-\frac{h}{2},y - \frac{\sqrt{3}}{2}h)}{8\sqrt{3}} \\
-c^2 \alpha_2  \frac{3 u(x,y) + u(x-h,y) + u(x-\frac{h}{2},y + \frac{\sqrt{3}}{2}h) + u(x+\frac{h}{2},y + \frac{\sqrt{3}}{2}h) + u(x+h,y) + u(x+\frac{h}{2},y - \frac{\sqrt{3}}{2}h) +u(x-\frac{h}{2},y - \frac{\sqrt{3}}{2}h)}{6\sqrt{3}} =
 \\
-\frac{189 h^4 u^{(0,6)}(x,y)+945 h^4 u^{(2,4)}(x,y)+315 h^4 u^{(4,2)}(x,y)+231 h^4 u^{(6,0)}(x,y)+7560 h^2 u^{(0,4)}(x,y)+15120 h^2 u^{(2,2)}(x,y)+7560 h^2 u^{(4,0)}(x,y)+120960 u^{(0,2)}(x,y)+120960 u^{(2,0)}(x,y)}{80640 \sqrt{3}}
\\ 
-\frac{c^2 \left(h^2 \left(7560 h^2 u^{(0,4)}(x,y)+15120 h^2 u^{(2,2)}(x,y)+7560 h^2 u^{(4,0)}(x,y)+120960 u^{(0,2)}(x,y)+120960 u^{(2,0)}(x,y)\right)+967680 u(x,y)\right)}{645120 \sqrt{3}} \\
-\frac{\text{$\alpha $2} c^2 \left(h^2 \left(7560 h^2 u^{(0,4)}(x,y)+15120 h^2 u^{(2,2)}(x,y)+7560 h^2 u^{(4,0)}(x,y)+120960 u^{(0,2)}(x,y)+120960 u^{(2,0)}(x,y)\right)+725760 u(x,y)\right)}{483840 \sqrt{3}} 
  + \mathcal{O}(h^5).
\end{array}
 \label{eqn:truncationerror2d}
\end{eqnarray}
Since 
\begin{equation*}
 -\frac{120960 u^{(0,2)}(x,y)+120960 u^{(2,0)}(x,y)}{80640 \sqrt{3}}- \frac{967680 u(x,y)}{645120 \sqrt{3}} = 0,
\end{equation*}
(\ref{eqn:truncationerror2d}) becomes
\begin{eqnarray}
 \begin{array}{ll}
\tau(x,y)  =
-\frac{189 h^4 u^{(0,6)}(x,y)+945 h^4 u^{(2,4)}(x,y)+315 h^4 u^{(4,2)}(x,y)+231 h^4 u^{(6,0)}(x,y)+7560 h^2 u^{(0,4)}(x,y)+15120 h^2 u^{(2,2)}(x,y)+7560 h^2 u^{(4,0)}(x,y)}{80640 \sqrt{3}}
\\ 
-\frac{c^2 \left(h^2 \left(7560 h^2 u^{(0,4)}(x,y)+15120 h^2 u^{(2,2)}(x,y)+7560 h^2 u^{(4,0)}(x,y)+120960 u^{(0,2)}(x,y)+120960 u^{(2,0)}(x,y)\right)\right)}{645120 \sqrt{3}} \\
-\frac{\text{$\alpha $2} c^2 \left(h^2 \left(7560 h^2 u^{(0,4)}(x,y)+15120 h^2 u^{(2,2)}(x,y)+7560 h^2 u^{(4,0)}(x,y)+120960 u^{(0,2)}(x,y)+120960 u^{(2,0)}(x,y)\right)+725760 u(x,y)\right)}{483840 \sqrt{3}} 
  + \mathcal{O}(h^5).
\end{array}
 \label{eqn:truncationerror2d2}
\end{eqnarray}
Substituting the exact solution $u(x,y) = \sin(c\cos(\theta)x + c\sin(\theta)y)$ (note that more general solutions can be chosen) into (\ref{eqn:truncationerror2d2})  gives
\begin{eqnarray}
 \begin{array}{ll}
\tau(x,y)  =
 -\frac{1}{322560 \sqrt{3}} (c^2 (35 c^5 h^5 \cos (c x \cos (\theta )+c y \sin (\theta )+\theta )+21 c^5 h^5 \cos (c x \cos (\theta )+c y \sin (\theta )+3 \theta )+7 c^5 h^5 \cos (c x \cos (\theta )+c y \sin (\theta )+5 \theta )+ \\ 
 c^5 h^5 \cos (c x \cos (\theta )+c y \sin (\theta )+7 \theta )+35 c^5 h^5 \cos (-c x \cos (\theta )-c y \sin (\theta )+\theta )+21 c^5 h^5 \cos (-c x \cos (\theta )-c y \sin (\theta )+3 \theta )+7 c^5 h^5 \cos (-c x \cos (\theta )-c y \sin (\theta )+5 \theta )+ \\
 c^5 h^5 \cos (-c x \cos (\theta )-c y \sin (\theta )+7 \theta )+5040 \text{$\alpha $2} c^4 h^4 \sin (c x \cos (\theta )+c y \sin (\theta ))+2940 c^4 h^4 \sin (c x \cos (\theta )+c y \sin (\theta ))-42 c^4 h^4 \sin (c x \cos (\theta )+c y \sin (\theta )+6 \theta )+ \\ 
 42 c^4 h^4 \sin (-c x \cos (\theta )-c y \sin (\theta )+6 \theta )-80640 \text{$\alpha $2} c^2 h^2 \sin (c x \cos (\theta )+c y \sin (\theta ))-30240 c^2 h^2 \sin (c x \cos (\theta )+c y \sin (\theta ))+483840 \text{$\alpha $2} \sin (c x \cos (\theta )+c y \sin (\theta)) 
  + \mathcal{O}(h^5).
\end{array}
 \label{eqn:truncationerror2d3}
\end{eqnarray}
The last 6 terms in (\ref{eqn:truncationerror2d3}) can be written in the form 
\begin{equation}
 \sin (c x \cos (\theta )+c y \sin (\theta)) \left(c^6h^4 \underbrace{\frac{-2940 + \frac{80640\alpha2}{c^2} +84 \cos(6\theta)}{322560 \sqrt{3}}}_{C_2}  + c^4h^2 \underbrace{\frac{30240 - \frac{483840 \alpha2}{c^2h^2}}{322560 \sqrt{3}}}_{C_1} \right). 
 \label{eqn:twocoefficient}
\end{equation}

\end{sidewaysfigure}

\begin{figure}[H]
\center
 \includegraphics[width=8cm]{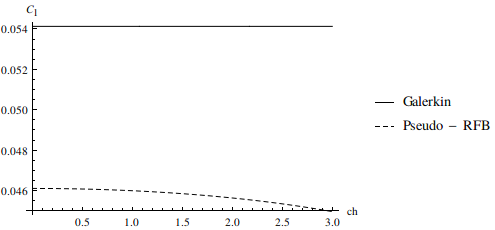} 
  \includegraphics[width=8cm]{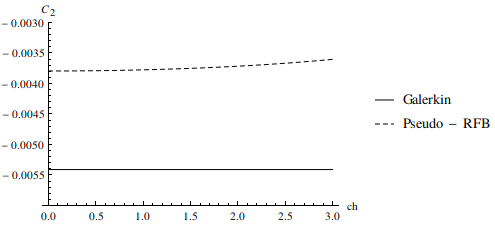} 
\caption{Comparison of the linear Galekin method and the pseudo-RFB method for the coefficients $C_1$ (left) and $C_2$ (right).}
\label{fig:coefficient1}
\end{figure}

Figure \ref{fig:coefficient1} shows the graph of $C_1$ and $C_2$  for $\alpha_2=0$ (standard Galerkin) and  for $\alpha_2$  in (\ref{eqn:alfagalerkin}) (pseudo-RFB) . We set $\theta=\pi/3$ to plot the graph of $C_2$.   The  slight decreases in $C_1$ and $C_2$ in magnitude  for $0<ch<3$,   explains why the RFB method is not effective in 2D.       
One way to  improve the accuracy of the  bubble  approach in 2D is to modify the right-hand side of the bubble equations in (\ref{eqn:2dbubble3}) by multiplying with a constant, say $\mu$. Then, $\alpha_2$ becomes 
\begin{equation*} 
 \alpha_2 = \frac{2\mu c^2h^2}{3(72-c^2h^2)}. 
\end{equation*}
After this modification, the bubble functions are no more residual-free. We call these modified functions as \textit{adaptive bubble} functions. We call the piecewise-defined linear approximations to these adaptive bubble functions as \textit{pseudo-adaptive} bubble functions. 

\begin{figure}[H] 
\center
 \includegraphics[width=8cm]{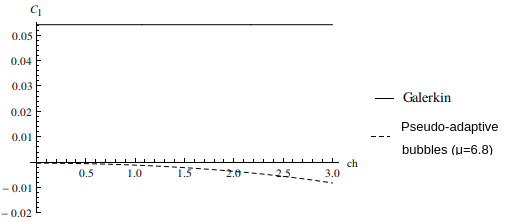} 
  \includegraphics[width=8cm]{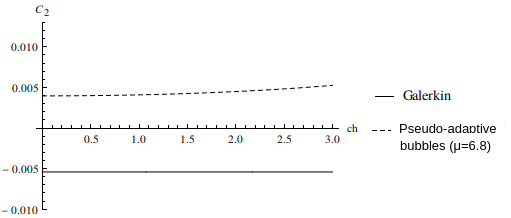} 
\caption{Comparison of the linear Galekin method and the pseudo-adaptive bubbles  method ($\mu=6.8$) for the coefficients $C_1$ (left) and $C_2$ (right).}
\label{fig:coefficient2}
\end{figure}

\begin{figure}[H]
\center
 \includegraphics[width=8cm]{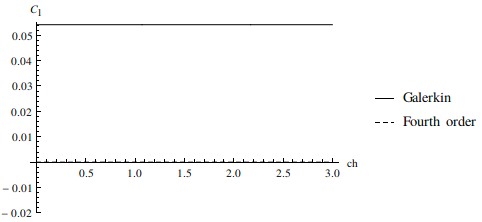} 
  \includegraphics[width=8cm]{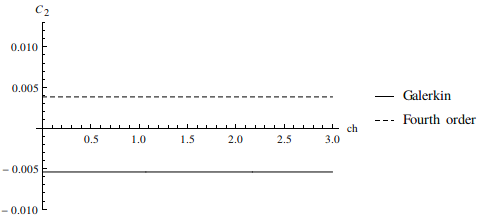} 
\caption{Comparison of the linear Galekin method and the fourth order accurate method for the coefficients $C_1$ (left) and $C_2$ (right).}
\label{fig:coefficient3}
\end{figure}

We give two examples here to validate the approach.  
Figure \ref{fig:coefficient2} shows the graph of $C_1$ and $C_2$ when $\mu=6.8$.  It is clear that $C_1$ is decreased in magnitude substantially. It is almost zero when $ch$ is close to zero. There is not much change in $C_2$ in magnitude.  
For the second example, consider  $\alpha_2=0.0625 c^2h^2$ which  makes $C_1$ zero for all values of $ch$. In this case, the finite difference scheme in  (\ref{eqn:2dfinitedifference}) is a fourth order scheme with seven points for plane waves. 
Figure \ref{fig:coefficient3} shows the graph of $C_1$ and $C_2$ when $\alpha_2=0.0625 c^2h^2$. While $C_1$ is zero for all values of $ch$, there is only a slight change in $C_2$ in magnitude. This fourth order accurate finite difference scheme can be easily applied  in triangular, trapezoidal  and polygonal domains. Our main aim in this article is to propose adaptive bubbles approximated by standard Galerkin method. However, the above two method will be used to  compare the success of the AB method.      

\section{Adaptive bubbles (AB) in 2D with triangular elements} \label{sec:section4}
We have shown using the pseudo-bubbles and the truncation error that the RFB method is not effective in 2D. However, it is possible to increase its accuracy with a simple modification, that is, multiplying the right-hand side of the bubble problems with a constant. We proposed two methods using this approach; a pseudo-adaptive bubbles method and a fourth order accurate finite difference scheme that uses seven points. However, our main aim is to obtain more accurate solutions by approximating the adaptive bubble functions with linear finite element method on a coarse mesh.         

\begin{figure}[H]
 \center
\includegraphics[width=8cm]{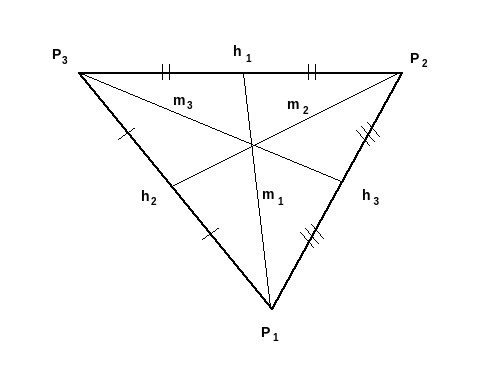}
\includegraphics[width=8cm]{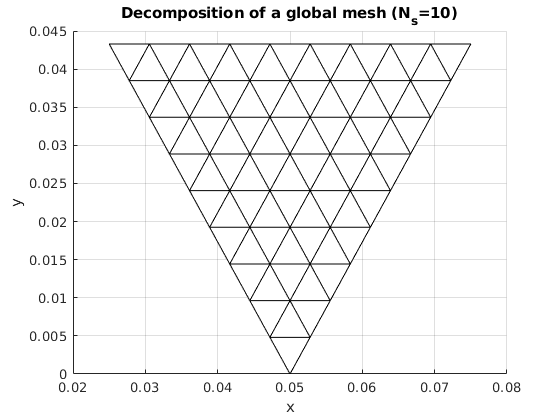}
\caption{A global mesh (left) and a decomposition of a global mesh with triangular elements when $N_s=10$ (right).}
\label{fig:triangle}
\end{figure}

We follow an empirical way to determine the optimal values of $\mu_i$ for varying  $cm_i$ where $m_i$ is the median of the global triangular element (see Figure \ref{fig:triangle}, (left)). This is actually a necessity because we have to use a coarse mesh for the sub-problems and   shapes of the bubble functions may change significantly when a small change occurs in the number of mesh used for the sub-problems.  
Considering the problem in (\ref{eqn:Helmholtz2ddirichlet}),
we report the optimal values of $\mu_i$ ($i=1,2,3$) for different values of $cm_i$ in Table \ref{table:optimal} for an equilateral triangular  element. The  optimality criterion in determining these values is minimization of the error in inifinity norm by doing many tests.  The values of $\mu_i$ are optimal for $\theta=0,\pi/3, 3\pi/2$.   We decompose each  global mesh into triangular elements by choosing $N_s$ uniformly distributed  nodes on all edges of a triangular element (see Figure \ref{fig:triangle}).      
We set $N_s=10$ when $cm_i \leq 2.577$ and set $N_s=15$ when $cm_i> 2.577$. While choosing $N_s=10$ amounts to   $28 \times 28$  matrices, choosing $N_s=15$ amounts to   $78 \times 78$ matrices on element level. It is possible to choose smaller values of $N_s$, especially for smaller $ch$, but one has to report more optimal values of $\mu$ in this case.     
\begin{table}[H]
\caption{Optimal values of $\mu_i$ for varying $cm_i$ ($i=1,2,3$) when $\theta=0,\pi/3, 3\pi/2$.}  
\begin{center}
\begin{tabular}{ll}
\begin{tabular}{|c|c|c|}
\hline
$cm_i$ & $\mu_i$ & $N_s$  \\
\hline
$\leq$ 0.57 & 5.4 & 10 \\ 
0.583 & 5.43 & 10 \\ 
0.644 & 5.45 & 10 \\ 
0.71 & 5.5 & 10 \\
0.7876 & 5.5 & 10  \\
0.859 & 5.51  & 10 \\
0.893 & 5.58  & 10 \\
0.930 & 5.6 & 10 \\
1.002 & 5.65 & 10 \\
1.074 & 5.7 & 10 \\
1.145 & 5.75 & 10 \\
1.217 & 5.8 & 10 \\
1.288 & 5.87 & 10 \\
1.360 & 5.95 & 10 \\
1.431 & 6.05 & 10 \\
1.503 & 6.1 & 10 \\
1.575 & 6.17 & 10 \\
1.646 & 6.3 & 10 \\
1.718 & 6.4 & 10 \\
1.789 & 6.5 & 10 \\
1.861 & 6.6 & 10 \\
\hline
\end{tabular}
&
\begin{tabular}{|c|c|c|}
\hline
$cm_i$ & $\mu_i$ & $N_s$  \\
\hline
1.933 & 6.75 & 10 \\
2.004 & 6.9 & 10 \\
2.076 & 7.05 & 10 \\
2.147 & 7.2 & 10 \\
2.219 & 7.35 & 10 \\
2.291 & 7.52 & 10 \\
2.362 & 7.7 & 10 \\
2.434 & 7.9 & 10 \\
2.505 & 8.05 & 10 \\
2.577 & 8.3 & 10 \\
2.577 & 7.8 & 15 \\
2.649 & 7.95 & 15 \\
2.72  & 8.1  & 15 \\
2.75  & 8.21  & 15 \\
2.79  & 8.25 & 15 \\
2.863 & 8.4  & 15\\
2.93  & 8.5  & 15 \\
3.007 & 8.55 & 15 \\
3.078 & 8.6  & 15 \\
3.15  & 8.65 & 15 \\
\hline
\end{tabular}
\end{tabular}
\end{center}
\label{table:optimal}
\end{table}
When $cm_i$ is between any of the  	successive	two  values in Table \ref{table:optimal}, we use linear  interpolation to get $\mu_i$. When shape of the global mesh changes,   the bubble functions behave differently, and hence it becomes more difficult to find the optimal values of $\mu_i$ for each bubble functions. We therefore expect deterioration of the AB  method when nonuniform mesh is used, especially for large wave numbers. Rectangular elements require solving 5 different bubble problems on each element and hence it becomes more complicated to determine the optimal values of $\mu_i$ on nonuniform mesh. Hence, we expect the triangular elements to be more efficient than the rectangular elements on nonuniform mesh.

\remark \textit{The optimal values in Table \ref{table:optimal} were determined when $\theta=0$ and hence they are true values for $\theta=\pi/3, 2\pi/3$. The optimal values can be find for other values of $\theta$. We will show by numerical test that the values of $\mu_i$ in Table  \ref{table:optimal}  can be used in any direction when $ch<1$. It is possible to obtain good approximations up to $ch=2$ in any direction when $c<200$. }

\remark \textit{When $cm_i<0.57$ (i,e. $ch<0.65$) the optimal values in any direction are same, that is, $\mu=5.4$. In simulations, 10 nodes per wave are generally used which corresponds to $ch\approx 0.625$. In this regime, there is only one parameter that we must use, that is, $\mu=5.4$. Since it works for every direction, we expect the AB method works efficiently when the solution is not a plane wave or an unstructured mesh is used.}

\begin{figure}[H]
\center
 \includegraphics[width=8cm]{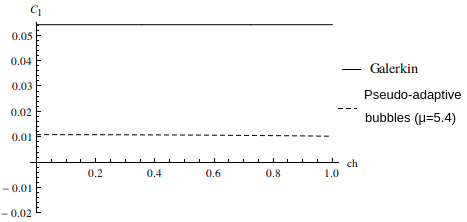} 
  \includegraphics[width=8cm]{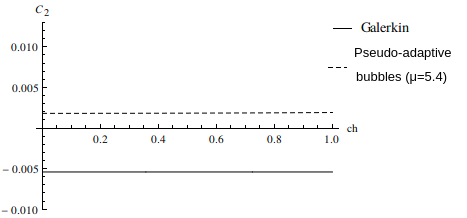} 
\caption{Comparison of the linear Galekin method and the pseudo-adaptive bubbles  method ($\mu=5.4$) for the coefficients $C_1$ (left) and $C_2$ (right).}
\label{fig:coefficient5}
\end{figure}

To verify the optimal values in Table \ref{table:optimal}, we use the pseudo-adaptive bubbles. For example, for $\mu=5.4$ (when $ch<0.65$), $\alpha_2=10.8c^2h^2/(3(72-c^2h^2))$.  Graphs of the coefficients $C_1$ and $C_2$ are provided in Figure \ref{fig:coefficient5}.  It is obvious that both $C_1$ and $C_2$ are decreased in magnitude which is a verification that  the AB method can  mitigate the pollution effect substantially.     

\section{Numerical experiments}
\label{sec:section5}

In this section, we provide numerical tests to  asses the success of the AB method. We compare the AB method with   the pseudo-adaptive bubbles ($\mu=6.8$) (PAB),   RFB,  fourth order  and  standard linear  Galerkin methods.   We use standard linear Galekin method to approximate the bubble functions for the AB and RFB methods.

\subsection{Numerical test 1}

\begin{figure}[H] 
\center
 \includegraphics[width=8cm]{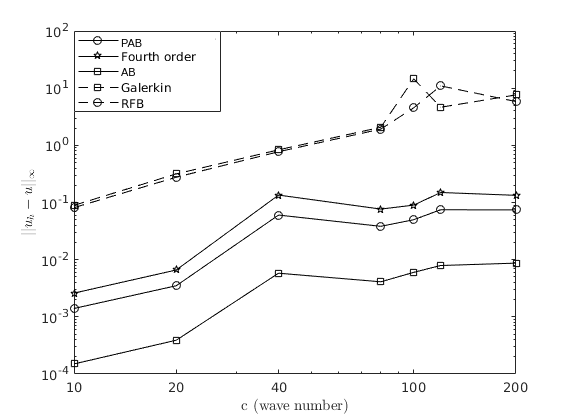} 
  \includegraphics[width=8cm]{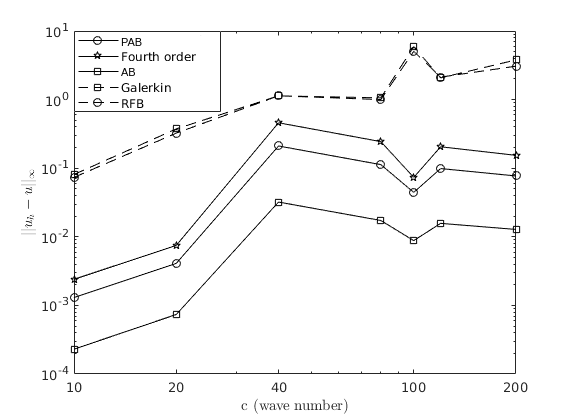} 
\caption{Comparison of the methods when $ch=0.625$ for $\theta=\pi/4$ (left) and $\theta=0$ (right). }
\label{fig:errorinf625}
\end{figure}

\begin{figure}[H]
\center
 \includegraphics[width=8cm]{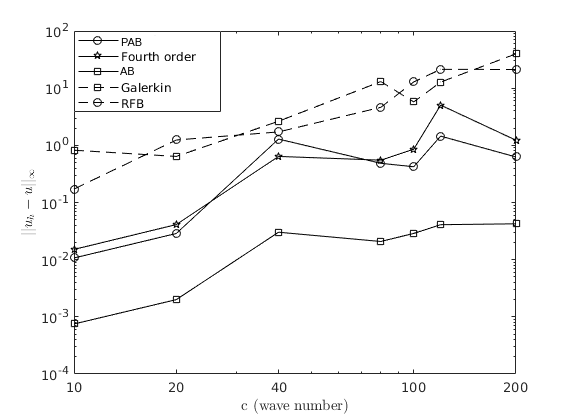} 
  \includegraphics[width=8cm]{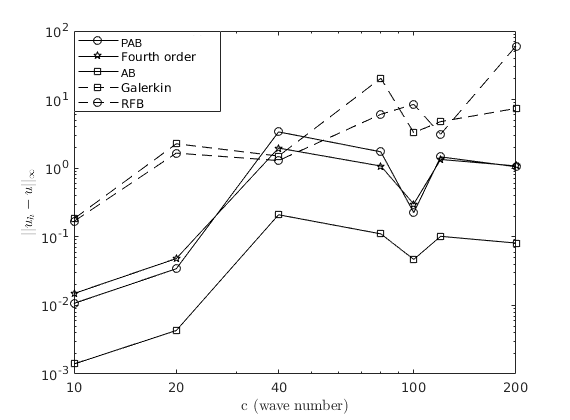} 
\caption{Comparison of the methods when $ch=1$ for $\theta=\pi/4$ (left) and $\theta=0$ (right). }
\label{fig:errorinf1}
\end{figure}

\begin{figure}[H]
\center
 \includegraphics[width=8cm]{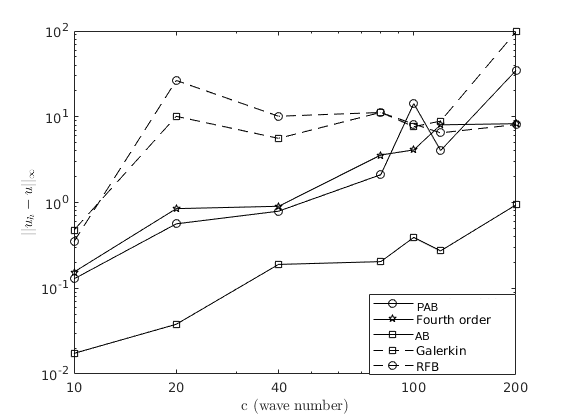} 
  \includegraphics[width=8cm]{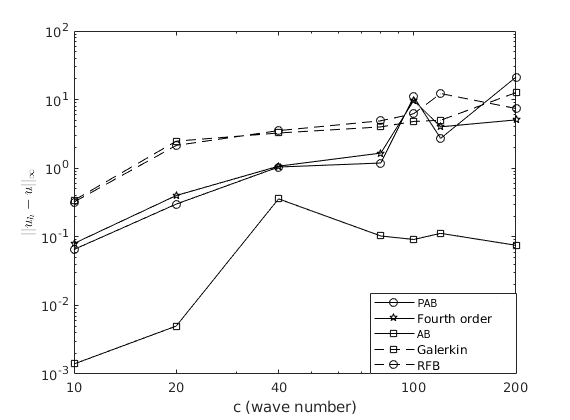} 
\caption{Comparison of the methods when $ch=1.75$ for $\theta=\pi/4$ (left) and $\theta=0$ (right). }
\label{fig:errorinf175}
\end{figure}

We consider the Helmholtz problem in (\ref{eqn:Helmholtz2ddirichlet}). Equilateral triangular elements are used to decompose the domain. We consider the cases  $ch=0.625,1,1.75$  when $\theta=0,\pi/4$ for increasing wave number to compare the methods in mitigating the pollution effect.   
Figure \ref{fig:errorinf625}, \ref{fig:errorinf1} and \ref{fig:errorinf175} show the log-log plots of the error in infinity norm for $ch=0.625$, $ch=1$ and $ch=1.75$, respectively. It is obvious that the AB method is better by far. The RFB method has  very small contribution in stabilization of the standard Galerkin method. The PAB method   outperforms the fourth order scheme.   
Moreover, the pollution error for the pseudo-bubbles method and the fourth order scheme is not negligible, particularly when $ch=1,1.75$. 

\begin{figure}[H]
\center
 \includegraphics[width=8cm]{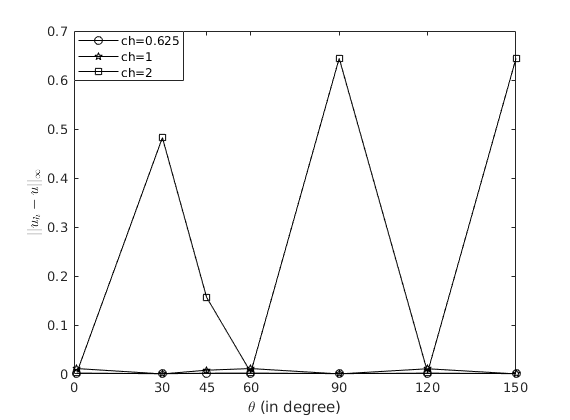} 
\caption{Error versus $\theta$ for the AB method when $c=50$. }
\label{fig:errorinfc50}
\end{figure}

Furthermore, we report errors for the M-RFB method in infinity norm for  $c=50$ and varying $\theta$ and $ch$ in Figure \ref{fig:errorinfc50}. It is clear that the direction of the plane waves has no importance in the error for $ch\leq 1$. One of the important observation is that the errors are almost same for $\theta=0,\pi/3, 2\pi/3$. A reasonable explanation for this is that $\cos(6x)$ appears as coefficient of $c^6h^4$ in the truncation error in  $C_2$ in (\ref{eqn:twocoefficient}).  This is directly related to the topology of the mesh. We can not expect the same behavior for rectangular elements.    

\subsection{Numerical test 2: Neumann boundary condition}
\begin{figure}[]
 \center
\includegraphics[width=6.5cm]{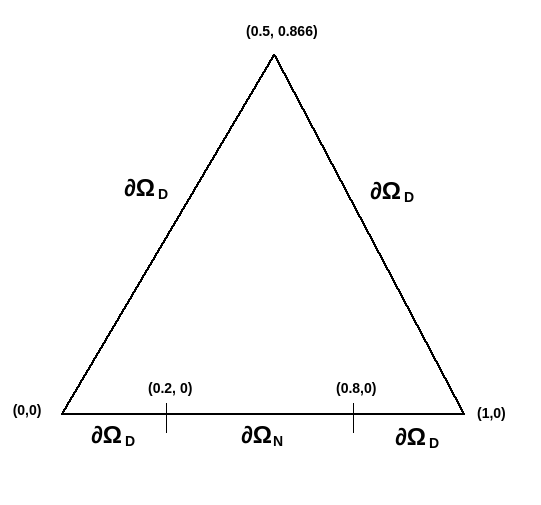}
\includegraphics[width=6cm, height=6cm]{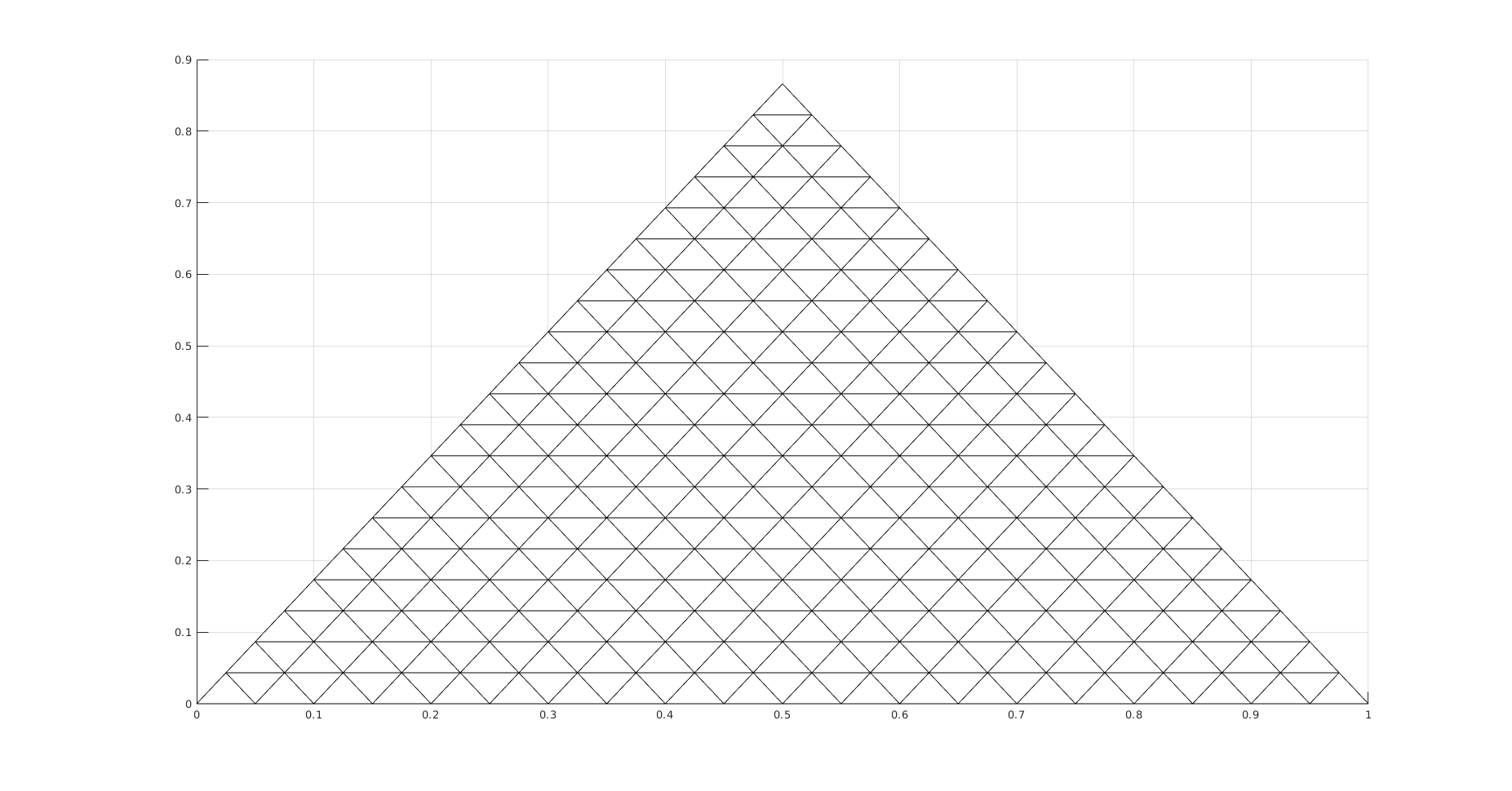}
\caption{Configuration of the domain in 2D (left) and its triangulation with equilateral triangular shaped elements (right).}
\label{fig:2ddomain}
\end{figure}

We consider the following  Helmholtz problem where homogenous Neumann boundary condition is imposed on a part of the boundary. 
\begin{eqnarray}
\left\lbrace
 \begin{array}{ll}
 -\Delta u -c^2u = 0, \quad \text{in} \quad \Omega, \\
  u(x,y) = \sin(cx),\quad \text{on} \quad \partial \Omega_D,\\
  \frac{\partial u}{\partial \textbf{n}} = 0,\quad \text{on} \quad \partial \Omega_N,
\end{array} \right.
 \label{eqn:Helmholtz2dint}
\end{eqnarray}
where $\Omega$, $\partial \Omega_D$ and $\partial \Omega_N$ are depicted in Figure \ref{fig:2ddomain}. Exact solution of (\ref{eqn:Helmholtz2dint}) is $u(x,y) = \sin(cx)$. To see the matrix formulation of the RFB method with Neumann boundary condition, we refer to \cite{RFBhelmholtz}. We decompose the domain with 400 equilateral  triangular  elements (see Figure \ref{fig:2ddomain}). 

\begin{figure}[H]
\center
 \includegraphics[width=14cm, height=7cm]{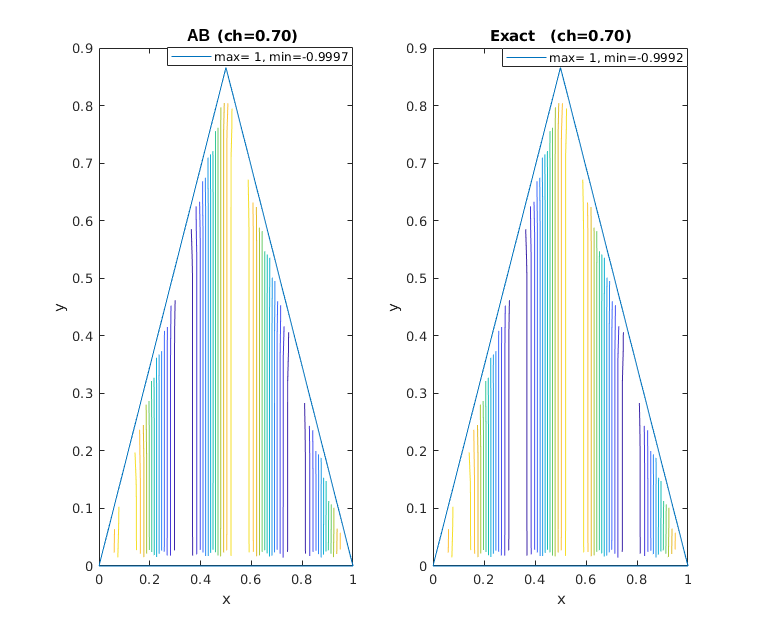}
 \includegraphics[width=14cm, height=7cm]{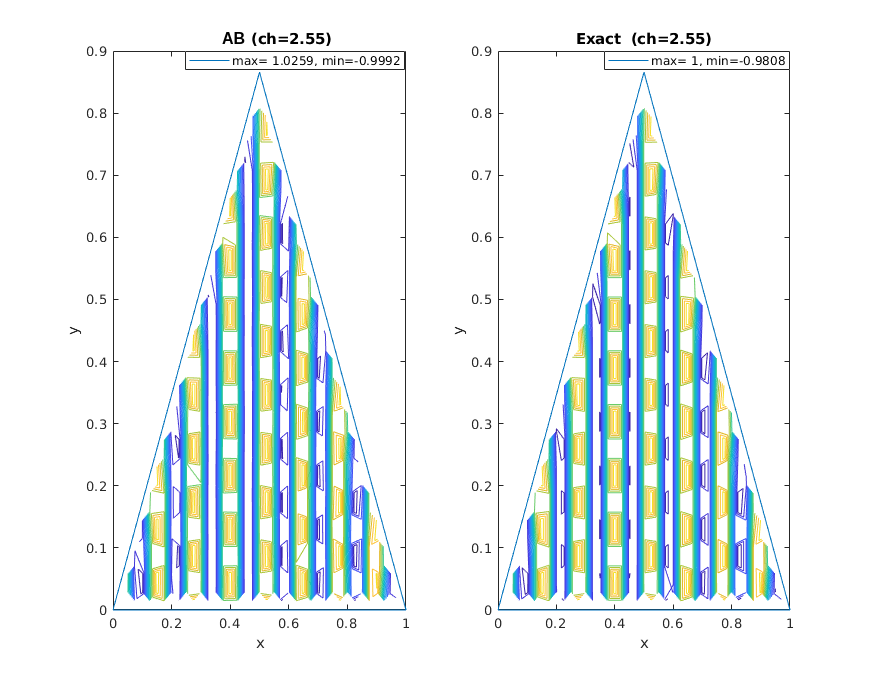}
  \includegraphics[width=14cm, height=7cm]{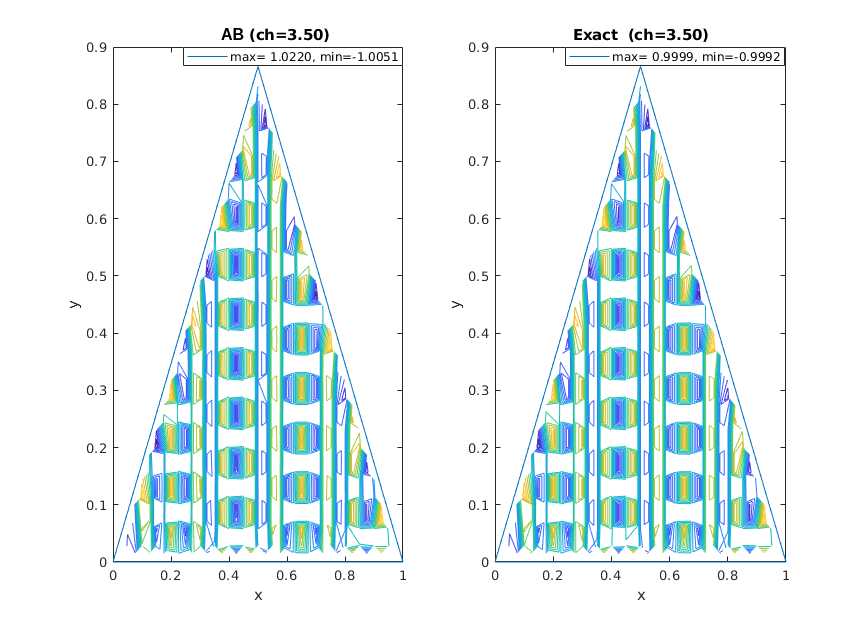}
\caption{Contour  plots of the approximate solutions obtained by the AB method and of the exact solutions for  $ch=0.70, 2.55, 3.50$.}
\label{fig:2dsolutionsMRFB2}
\end{figure}

\begin{figure}[H]
\center
 \includegraphics[width=18cm, height=7cm]{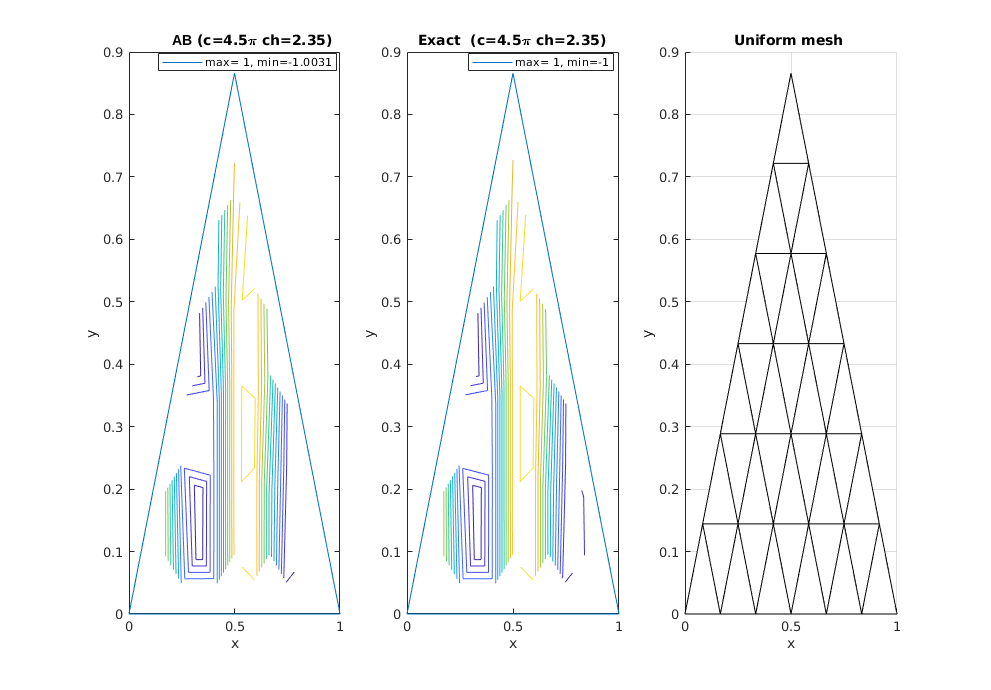}
 \includegraphics[width=18cm, height=7cm]{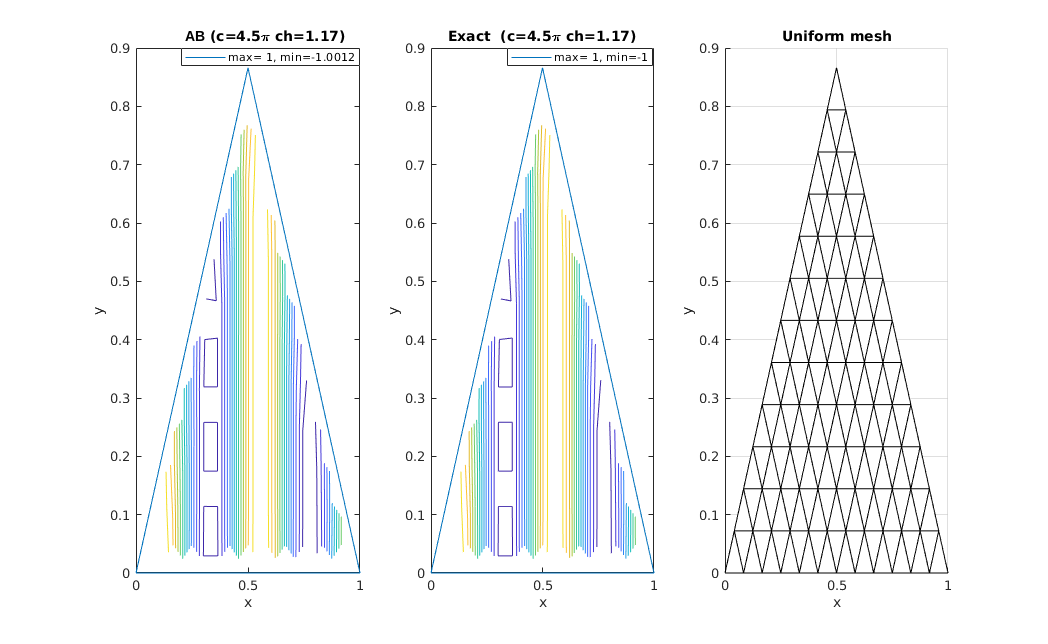}
  \includegraphics[width=18cm, height=7cm]{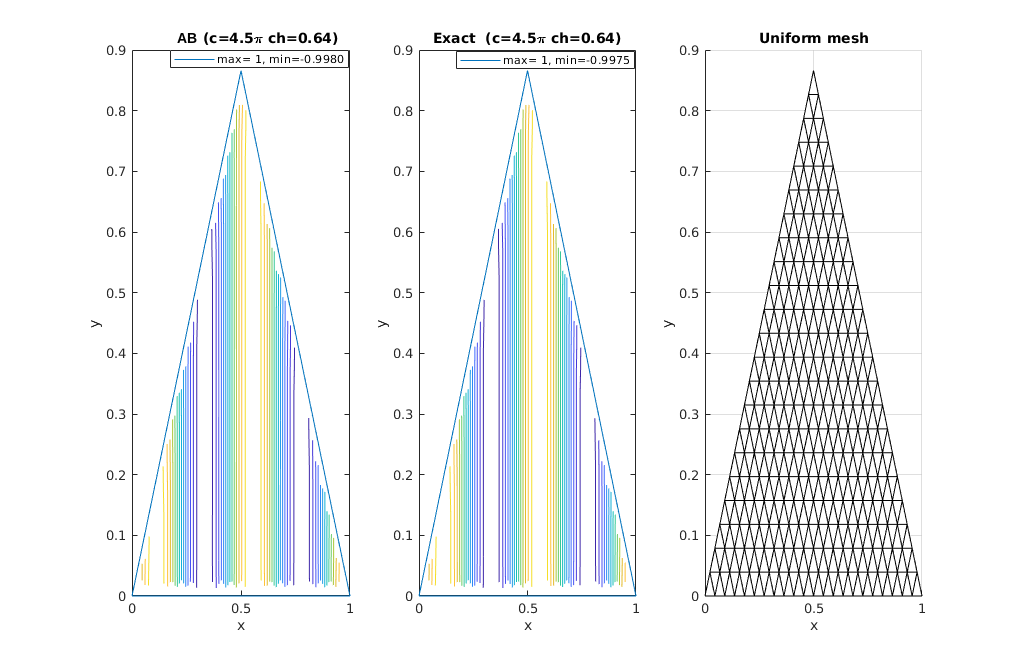}
\caption{Contour  plots of the approximate solutions obtained by the AB method on different  meshes and of the exact solutions on the same mesh  when  $c=4.5\pi$.}
\label{fig:2dsolutionsMRFB5}
\end{figure}

 We now test the AB method  with different values of $ch$. First, we report numerical results on a fixed uniform mesh (see Figure \ref{fig:2ddomain} (right)) for varying wave number $c$. Figures \ref{fig:2dsolutionsMRFB2} represents conotour plots of the approximate solutions and of the exact solutions for $ch=0.7,2.55, 3.5$. We also report  maximum and minumum values of the solutions on the graphs.  
Second, we report numerical results for fixed wave number $c$ on different meshes. Figure \ref{fig:2dsolutionsMRFB5} shows contour plots of the solutions and meshes used.

\begin{figure}[H]
\center
 \includegraphics[width=8cm, height=8cm]{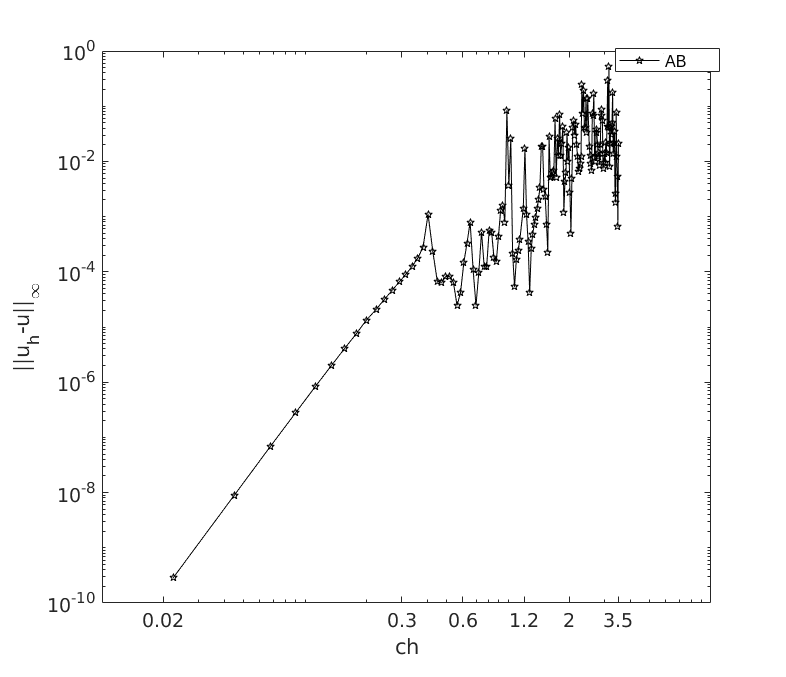}
\caption{Errors in infinity norm for the AB method up to $ch=3.5$.}
\label{fig:error}
\end{figure}

The results show that the AB method is very effective on uniform mesh up to $ch=3.50$. 
Finally, we report errors in inifinty norm in Figure \ref{fig:error} for the AB method up to $ch=3.5$ on a different mesh where 196 equilateral triangular elements are used. We calculated the error at many points. All the results above verify the robustness of the method in terms of the parameters proposed in Table \ref{table:optimal}.

% \begin{figure}[H]
% \center
%  \includegraphics[width=18cm, height=7cm]{2dMRFB070non.png}
%  \includegraphics[width=18cm, height=7cm]{2dMRFB086non.png}
%   \includegraphics[width=18cm, height=7cm]{2dMRFB125non.png}
% \caption{Contour  plots of the approximate solutions obtained by the M-RFB method on nonuniform mesh and of the exact solutions on the same mesh  when  $c=4.5\pi, 5.5\pi, 8\pi$.}
% \label{fig:2dsolutionsMRFBnon1}
% \end{figure}

\subsection{Numerical test 3:  Robin boundary condition and external source}

\begin{figure}[H]
\center
 \includegraphics[width=7cm, height=7cm]{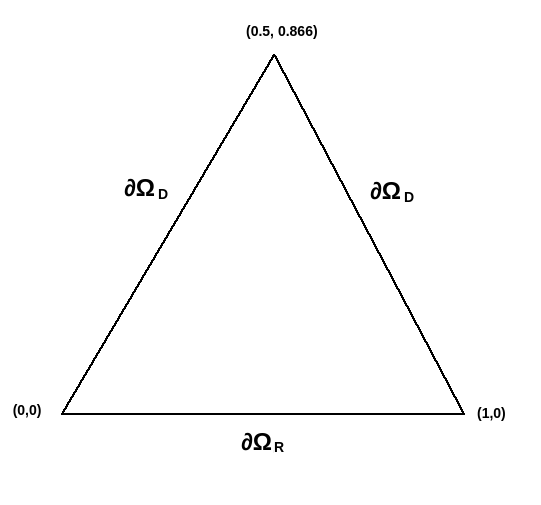}
  \includegraphics[width=7cm, height=7cm]{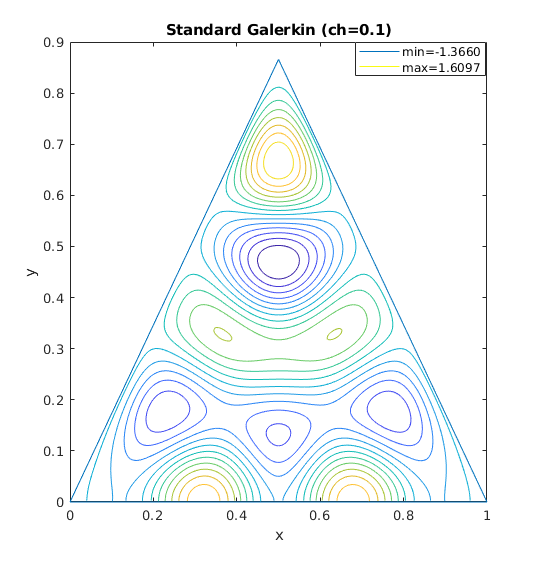}
\caption{Problem configuration (left) and a reference solution (right) obtained with standard Galekin method with 40000 uniform triangular elements for which $ch=0.1$.}
\label{fig:robinproblem}
\end{figure}

\begin{figure}[H]
\center
  \includegraphics[width=16cm, height=7cm]{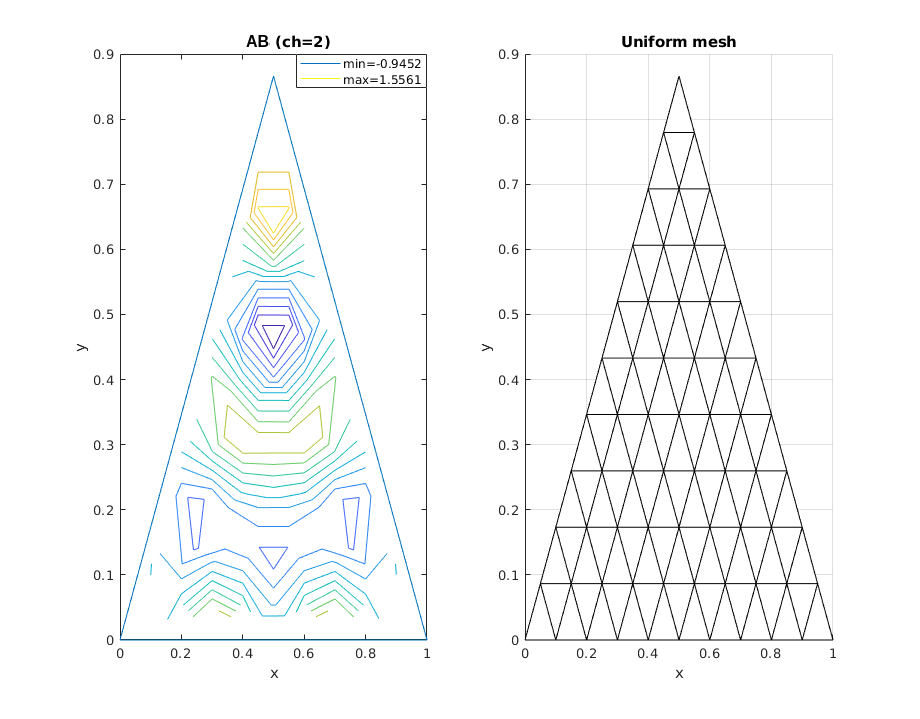}
 \includegraphics[width=16cm, height=7cm]{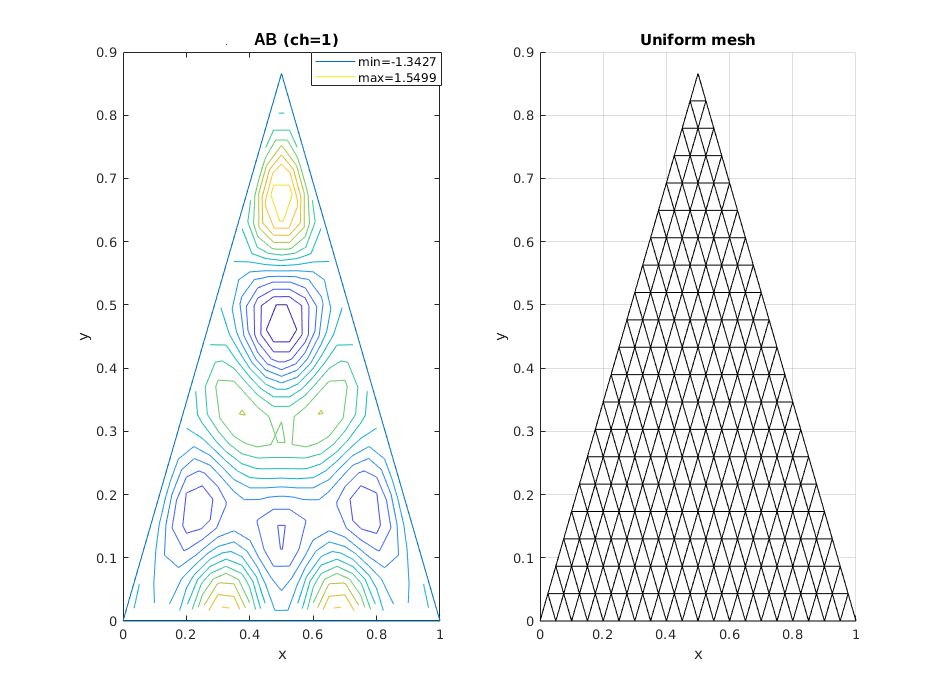}
 \includegraphics[width=16cm, height=7cm]{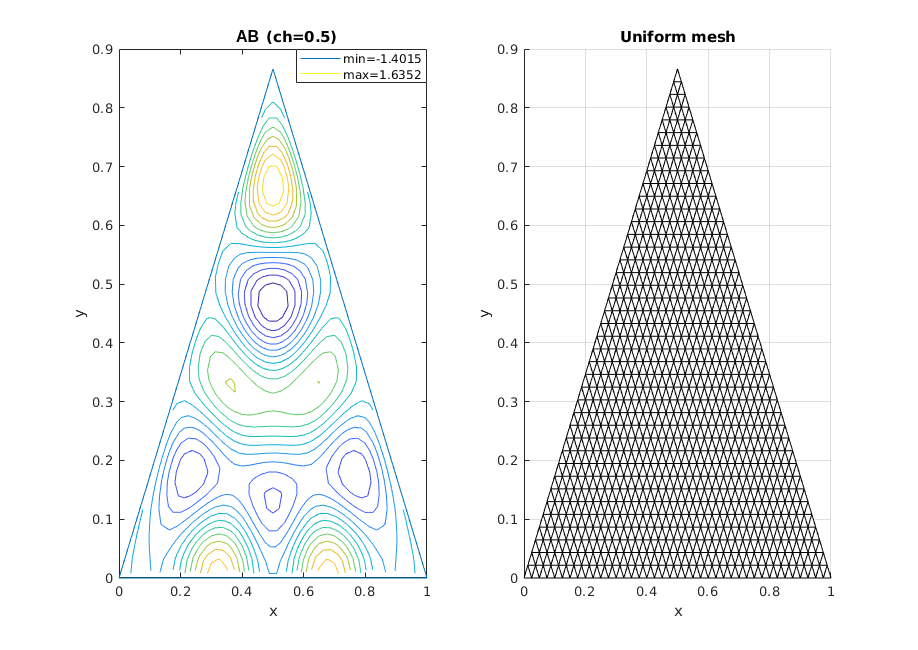}
\caption{Contour  plots of the approximate solutions obtained by the M-RFB method  when $c=20$ using different meshes  ($ch=0.5, 1,2$).}
\label{fig:robin}
\end{figure}

We test the AB method when  Robin boundary condition is imposed on a part of the boundary of the domain. We consider
\begin{eqnarray}
\left\lbrace
 \begin{array}{ll}
 -\Delta u -c^2u = \sin(x), \quad \text{in} \quad \Omega, \\
  u(x,y) = 0.1,\quad \text{on} \quad \partial \Omega_D,\\
  \frac{\partial u}{\partial \textbf{n}} = iu,\quad \text{on} \quad \partial \Omega_R,
\end{array} \right.
 \label{eqn:robin}
\end{eqnarray}
where $c=20$, and $\Omega$,  $\Omega_D$ and $\Omega_R$ are represnted in Figure \ref{fig:robinproblem} (left). As a refence solution, we get a solution using standard Galerkin  method on a fine mesh where  40000 uniform triangular elements are used for which $ch=0.1$.  Figure \ref{fig:robinproblem} (right)  shows the contour plot of the real part of the solution.  We show contour plots of the real part of the solutions  obtained by the AB method for $ch=0.5, 1, 2$ in Figure \ref{fig:robin}. Moreover, we show the corresponding meshes and report the maximun and minimum values of the approximate solutions. Results show that the AB method shows the characteristics of the reference solution for all cases. This is important in  application of the multigrid method as a solver.

\subsection{Numerical test 4: L-shaped domain and a different triangulation}

\begin{figure}[H]  
\center
  \includegraphics[width=7.5cm, height=7.5cm]{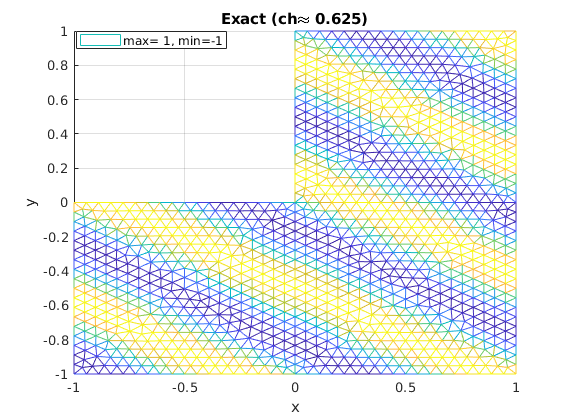}
   \includegraphics[width=7.5cm, height=7.5cm]{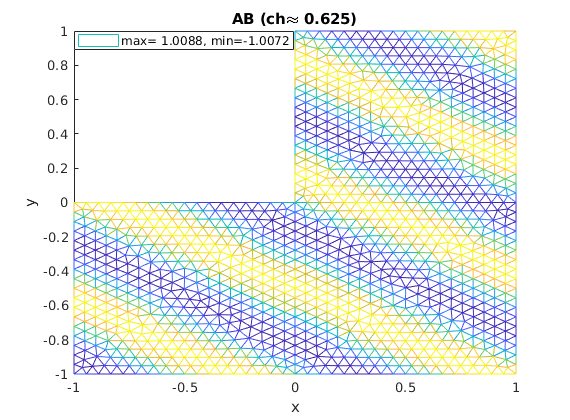}
   \includegraphics[width=7.5cm, height=7.5cm]{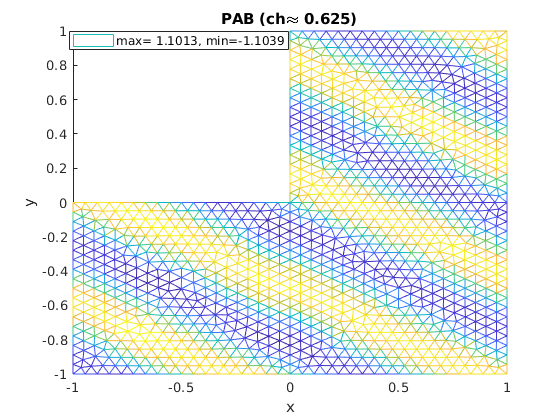}
    \includegraphics[width=7.5cm, height=7.5cm]{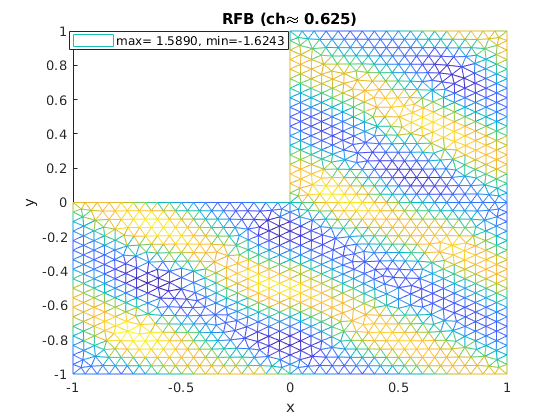}
\caption{Plots of the exact and approximate solutions obtained by the AB, PAB and RFB   methods  when $c=3.5\pi$.}
\label{fig:L-shaped}
\end{figure}

\begin{figure}[H]  
\center
  \includegraphics[width=7.5cm, height=7.5cm]{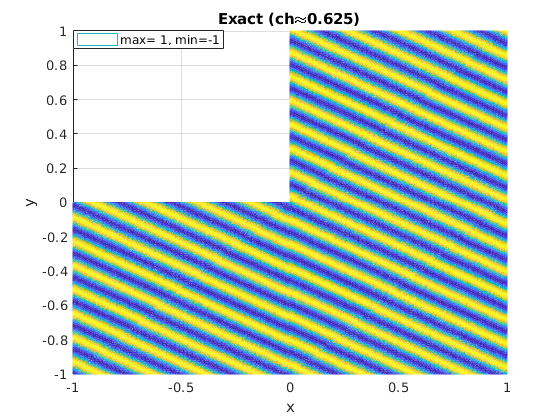}
   \includegraphics[width=7.5cm, height=7.5cm]{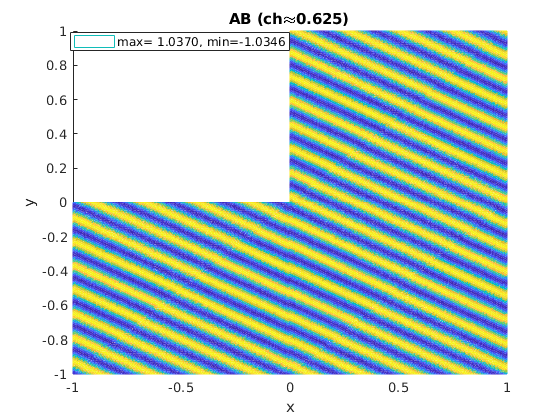}
   \includegraphics[width=7.5cm, height=7.5cm]{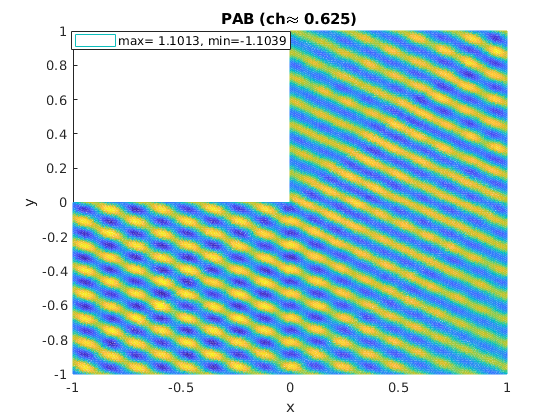}
    \includegraphics[width=7.5cm, height=7.5cm]{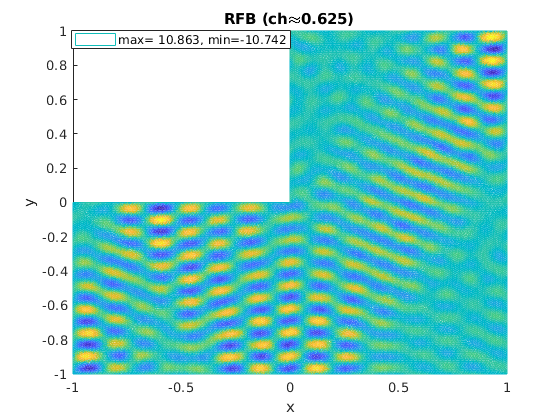}
\caption{Plots of the exact and approximate solutions obtained by the AB, PAB and RFB   methods  when $c=16.5\pi$..}
\label{fig:L-shaped2}
\end{figure}

\begin{figure}[H]  
\center
   \includegraphics[width=16cm, height=6.7cm]{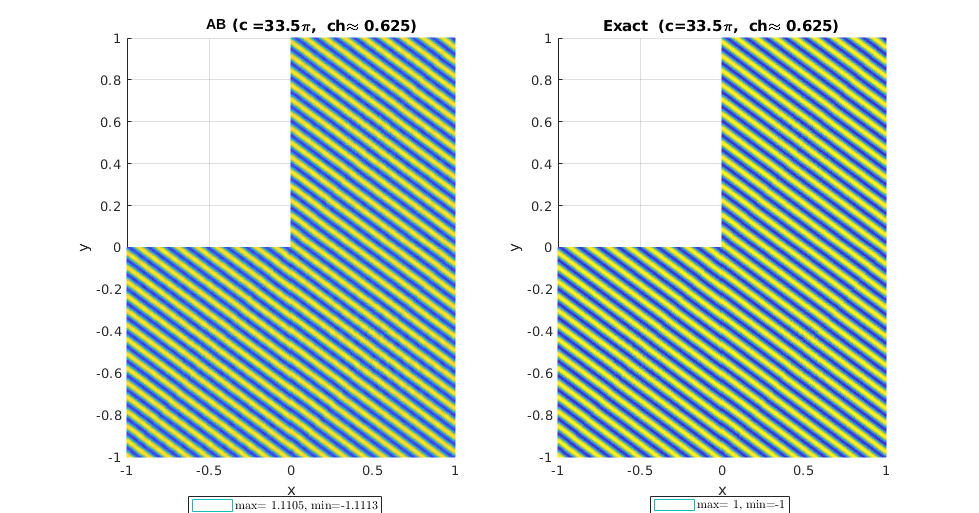}
\caption{Contour  plots of the approximate solutions obtained by the AB method and of the exact solutions on the same mesh when $c=33.5\pi$.}
\label{fig:L-shapedtrinagular}
\end{figure}

In this test problem, we change the domain and use a different triangulation. We use a L-shaped domain with the vertices $(-1,-1)$, $(-1,1)$, $(1,1)$, $(0,1)$, $(0,0)$ and $(-1,0)$. To decompose the domain, the following Matlab code  is used for which $ch\approx 0.625$.
\begin{eqnarray}
\left\lbrace
 \begin{array}{ll}
 model = createpde(1); \\
  geometryFromEdges(model,@lshapeg);\\
generateMesh(model,'GeometricOrder','linear','Hmax',0.625/c,'Hmin',0.625/c);
\end{array} \right.
 \label{}
\end{eqnarray}
The mesh  for the case $c=3.5\pi$ can be seen in Figure \ref{fig:L-shaped}. 
We consider the Dirichlet  problem in (\ref{eqn:Helmholtz2ddirichlet}) 
for $\theta=\pi/3$. Figure \ref{fig:L-shaped}, \ref{fig:L-shaped2} show the plots of the exact and approximate solutions obtained by the AB, PAB and RFB methods for $c=3.5\pi$ and $c=16.5\pi$, respectively. We also report the maximum and minimum values of the approximate solutions on the graphs. Results show that the AB method is better by far especially for larger wave number. Furthermore, we give the plots of the exact solution and approximate solution for AB given in Figure  \ref{fig:L-shapedtrinagular}. We did not report solutions for the PAB and RFB methods as their results are no more related to the exact solution.

\subsection{Numerical test 5: A circular complex domain}

\begin{figure}[H]  
\center
  \includegraphics[width=7.5cm, height=7.5cm]{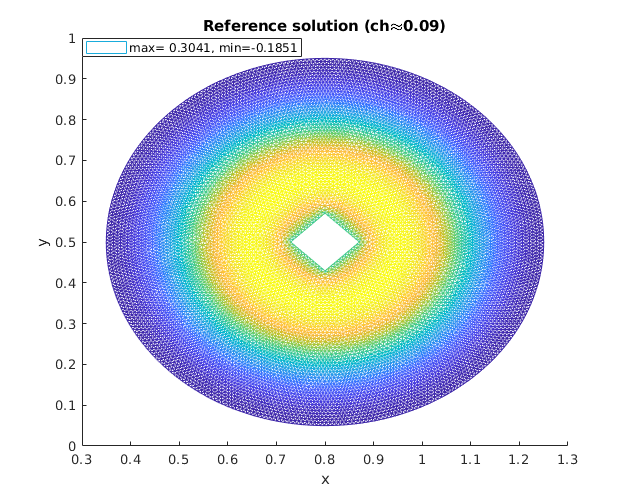}
   \includegraphics[width=7.5cm, height=7.5cm]{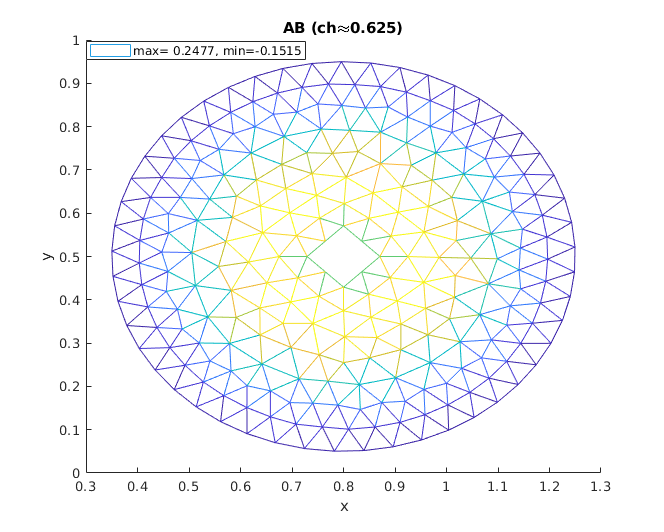}
   \includegraphics[width=7.5cm, height=7.5cm]{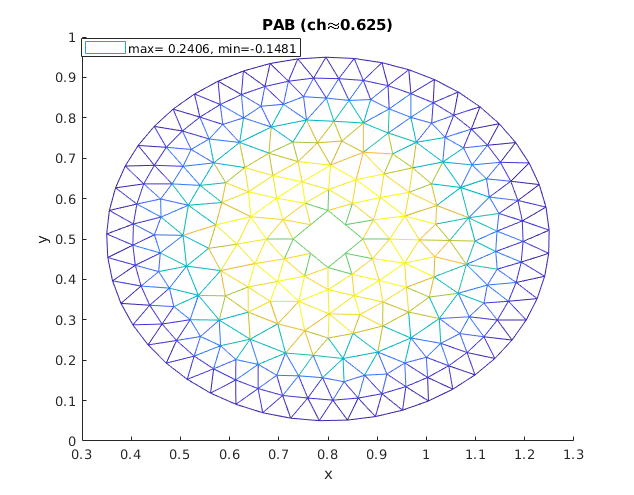}
    \includegraphics[width=7.5cm, height=7.5cm]{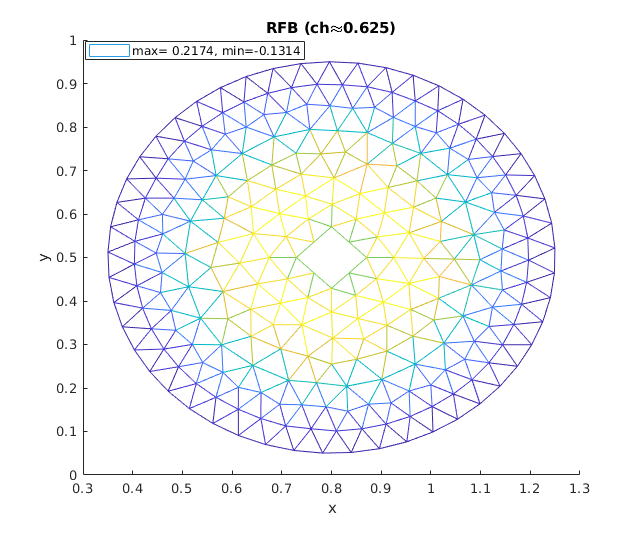}
\caption{Plots of  the reference solution and approximate solutions  obtained by the AB, PAB and RFB methods   when $c=3.5\pi$.}
\label{fig:circ-shaped}
\end{figure}

\begin{figure}[H]  
\center
  \includegraphics[width=7.5cm, height=7.5cm]{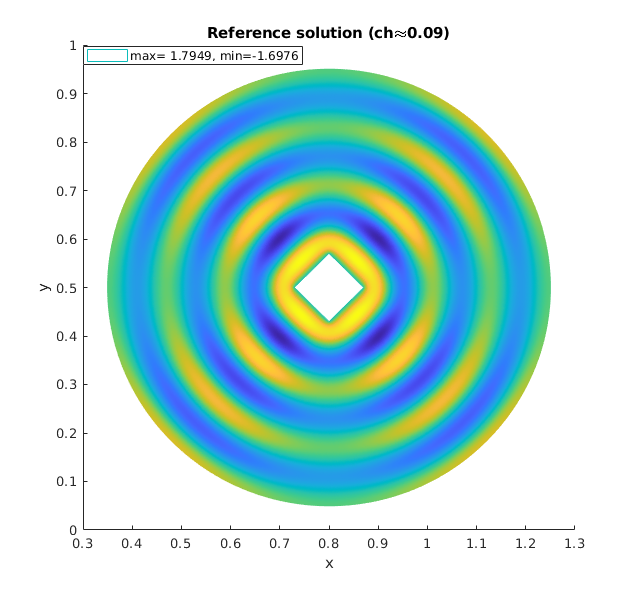}
   \includegraphics[width=7.5cm, height=7.5cm]{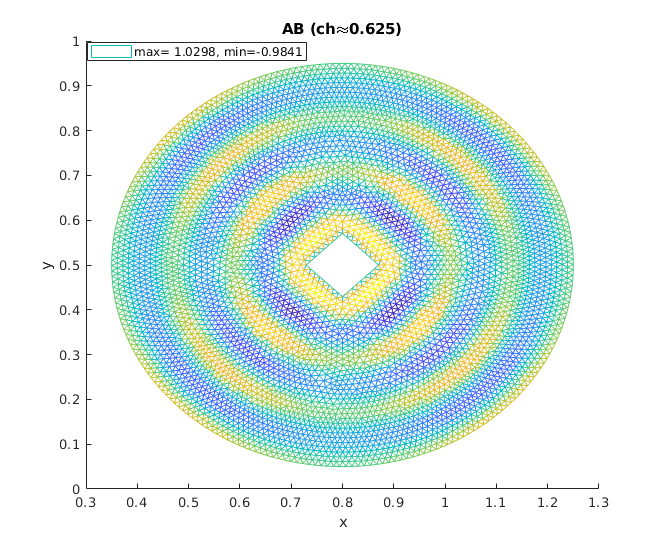}
   \includegraphics[width=7.5cm, height=7.5cm]{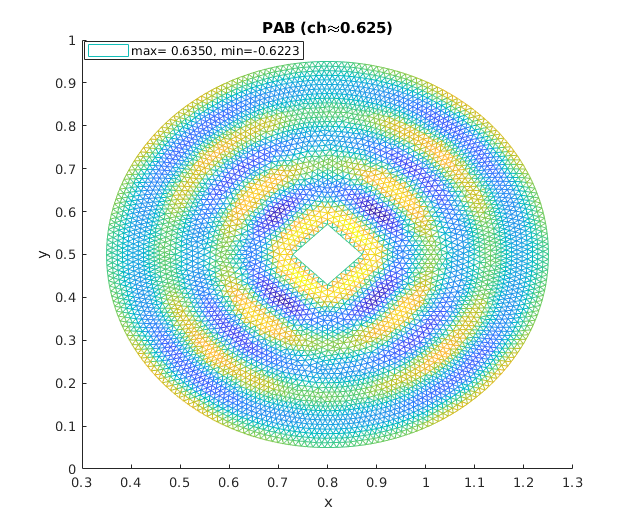}
    \includegraphics[width=7.5cm, height=7.5cm]{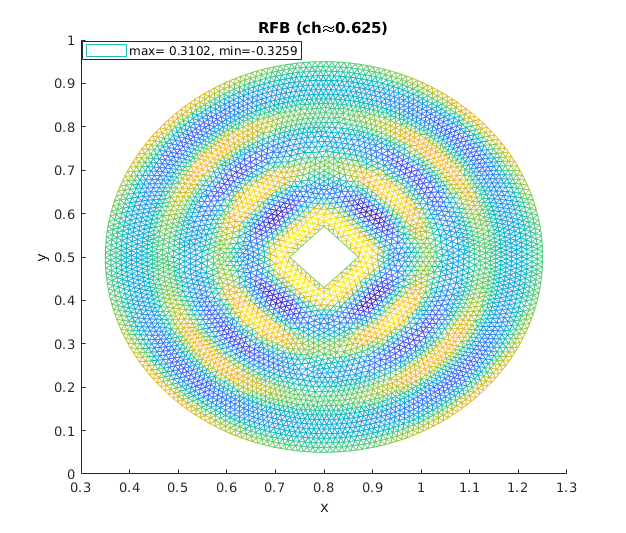}
\caption{Plots of  the reference solution and approximate solutions  obtained by the AB, PAB and RFB methods   when $c=16.5\pi$.}
\label{fig:circ-shaped2}
\end{figure}

In this test problem, we consider a complex domain which is obtained by the following Matlab code. 
\begin{eqnarray}
\left\lbrace
 \begin{array}{ll}
 model = createpde(1); \\
  geometryFromEdges(model,@scatterg);\\
generateMesh(model,'GeometricOrder','linear','Hmax',0.625/c,'Hmin',0.625/c);
\end{array} \right.
 \label{}
\end{eqnarray}
While homogenous Neumann boundary condition is imposed on the outer boundary of the domain,i.e, $\frac{\partial u}{\partial x} =0$, Dirichlet boundary condition is imposed on the inner boundary of the domain for which $u(x,y)=0.1$. 
The right hand side of the problem is set to zero.
The reference solution is obtained by standard Galerkin method on a fine mesh for which $ch\approx 0.09$. While Figure \ref{fig:circ-shaped} shows the plots of the reference solution and approximate solutions of the AB, PAB and RFB method for $c=3.5\pi$, Figure  \ref{fig:circ-shaped2} shows for $c=16.5\pi$. 
We see that the RFB method is worst in any case. Although the AB and PAB give similar results for smaller wave numbers, the AB method is better by far than the PAB  for large wave numbers. This and the previous tests show the success of the AB method on complex domains with unstructured meshes.

\section{Adaptive bubbles method with rectangular elements}  \label{sec:section6}

Although the AB method is very effective with triangular elements, in some domains,   rectangular elements may have some advantages such as in a rectangular region. For a rectangular element, there are five bubble equations to be solved.   
 \begin{eqnarray}
\left\lbrace
\begin{array}{ll}
  - \Delta \varphi_i  -c^2 \varphi_i =  \mu c^2 \psi_i \quad \text{in} \quad K, \quad (i=1,...,4)\\
 \varphi_i=0 \quad \text{on} \quad \partial K,
 \end{array}\right.
\label{eqn:2dbubble5} 
\end{eqnarray} 
and
 \begin{eqnarray}
\left\lbrace
\begin{array}{ll}
 - \Delta \varphi_f -c^2  \varphi_f = f \quad \text{in} \quad K,\\
  \varphi_f = 0 \quad \text{on} \quad \partial K,
 \end{array}\right.
\label{eqn:2dbubble6} 
\end{eqnarray} 
where $\psi_i$ $(i=1,...,4)$ are the bilinear basis functions of a rectangular element.  
In this case, the constant that we multiply the right hand side of the bubble equation (\ref{eqn:2dbubble5}) is fixed for each bubble equations. We considered the Dirichlet problem (\ref{eqn:Helmholtz2ddirichlet}) on unit square when $\theta=0$ to find the optimal values (in infinity norm). We report the optimal values of $\mu$ for a squared shaped element for varying $ch$ where $h=(h_1+h_2+h_3+h_4)/4$ and $h_i$ ($i=1,...,4$) are lengths of the edges of a rectangular element, in Table \ref{table:optimalrec}.

\remark \textit{Note that the values in Table \ref{table:optimalrec} are also optimal for $\theta=\pi$. It is possible to find the optimal values in any direction. However, triangular elements have some advantages.}
\begin{itemize}

 \item \textit{Rectangular elements use 9 points per degrees of freedom but triangular elements use  7 points per degrees of freedom}.
 
 \item \textit{While rectangular elements require  solving 5 different bubble equations, triangular elements requare 4. This makes rectangular elements less efficient when nonuniform mesh is used}.

 \item \textit{Triangular elements allow to work with larger $ch$}.

 \item \textit{Triangular elements are more efficient on unstructured meshes}. 

\end{itemize} 

\begin{table}[H] 
\caption{Optimal values of $\mu$ for rectangular elements for varying $ch$}  
\begin{center} 
\begin{tabular}{ll}
\begin{tabular}{|c|c|c|}
\hline
$ch$ & $\mu$ & $N_s$  \\
\hline
 $\leq$ 0.94 & 2.5  & 8 \\ 
1.02 & 2.6  & 8 \\ 
1.09 & 2.6  & 8 \\
1.17 & 2.6  & 8  \\
1.25 & 2.65 & 8 \\
1.33 & 2.65 & 8 \\
1.41 & 2.7  & 8 \\
1.49 & 2.7  & 8 \\
1.49 & 2.7  & 10 \\
1.57 & 2.72 & 10 \\
1.64 & 2.75 & 10 \\
\hline
\end{tabular}
&
\begin{tabular}{|c|c|c|}
\hline
$ch$ & $\mu$ & $N_s$  \\
\hline
1.72 & 2.8  & 10 \\
1.80 & 2.85 & 10 \\
1.88 & 2.88 & 10 \\
1.96 & 2.88 & 10 \\
2.04 & 2.98 & 10 \\
2.12 & 3.05 & 10 \\
2.19 & 3.09 & 10 \\
2.27 & 3.15 & 10 \\
2.35 & 3.2  & 10 \\
2.43 & 3.25 & 10 \\
2.51 & 3.29 & 10 \\
\hline
\end{tabular}
\end{tabular}
\end{center}
\label{table:optimalrec}
\end{table}

We are able to find the optimal values of $\mu$ up to $ch\approx2.5$. 
$N_s$ in Table \ref{table:optimalrec}, is the number of nodes on each edges of a recatngular element. $N_s=8$ for $ch\leq 1.49$ and  $N_s=10$  when $ch>1.49$. While $N_s=8$  amounts to solving  $36 \times 36$ linear  systems of equations, $N_s=10$  amounts to solving  $64 \times 64$ linear systems of equations on element level. When $ch$ is between any of the  	successive	two  values in Table \ref{table:optimalrec}, we use linear  interpolation to get $\mu$.
We provide one numerical test to show the performance of the AB.

\remark \textit{When $ch\leq 0.94$, the optimal values in any direction are same. Note that in simulations, 10 nodes per wave are generally used which corresponds to $ch\approx0.625$.}

\subsection{Test 1}
\begin{figure}[H]
 \center
\includegraphics[width=6.5cm]{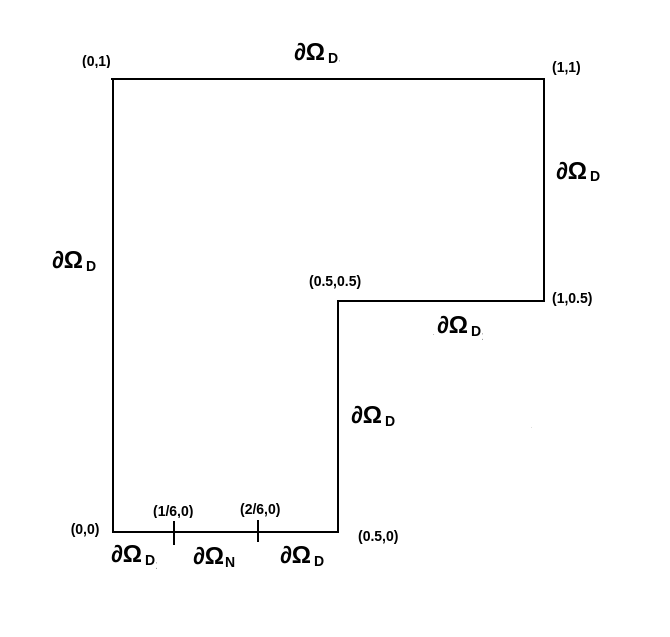}
\includegraphics[width=6cm, height=6cm]{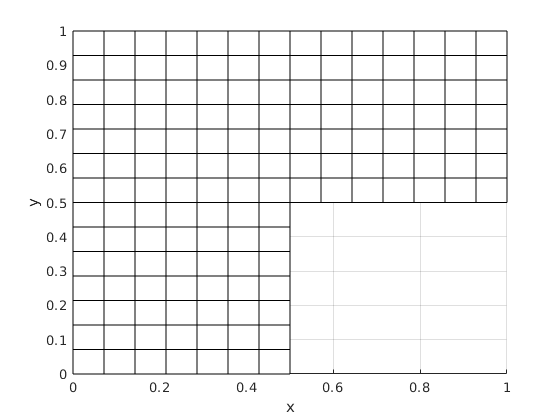}
\caption{Domain of the  problem (\ref{eqn:Helmholtz2lhaped}) (left) and its decomposition with square elements (right).}
\label{fig:2ddomainlshaped1}
\end{figure}

We consider the  following Helmholtz problem in 2D on an L-shaped domain (see Figure \ref{fig:2ddomainlshaped1} (left)).
\begin{eqnarray}
\left\lbrace
 \begin{array}{ll}
 -\Delta u -c^2u = 0, \quad \text{in} \quad \Omega, \\
  u(x,y) = \sin(cx),\quad \text{on} \quad \partial \Omega_D,\\
  \frac{\partial u}{\partial \textbf{n}} = 0,\quad \text{on} \quad \partial \Omega_N.
\end{array} \right.
 \label{eqn:Helmholtz2lhaped}
\end{eqnarray}
147 uniform square shaped  elements are used for the decomposition of the domain (see Figure \ref{fig:2ddomainlshaped1} (right)). Exact solution of this problem is $u(x,y)=\sin(cx)$.

\begin{figure}[H]
\center
  \includegraphics[width=16cm, height=7cm]{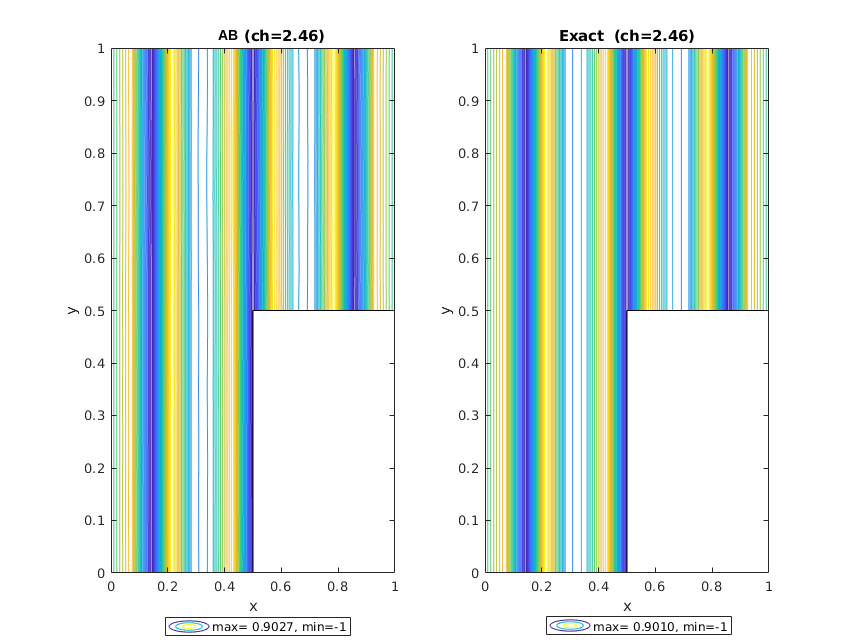}
\caption{Contour  plots of the approximate solutions obtained by the AB method  and of the exact solutions   when  $ch=2.46$.}
\label{fig:2dsolutionsMRFBrec1}
\end{figure}

We assess   the performance of the AB method by comparing with the exact solution.  Figure \ref{fig:2dsolutionsMRFBrec1} represents the contour plots of the approximate and of the exact solutions for $ch= 2.46$.  We also report the maximun and minimum values of the approximate and exact solutions.  Results show that the AB  method is very effective up to $ch=2.5$ on uniform mesh.

\section{Conclusion} \label{sec:section7}

In this article, we proposed   an adaptive bubble approch for the Helmhotz equation in 2D.  The RFB method requires obtaining the bubble functions which is generally as difficult as solving the orignal problem. We showed that this is not the case  for the Helmholtz problem.   The standard Galerkin finite element method can be used as a solver to obtain approximations to the  bubble functions. In other words, the bubbles method does not depend on another stabilized method when applied to the Helmholtz problem. We showed that the contribution of the RFB method in stabilization of the standard Galerkin method is very poor in 2D. We modified the RFB method by multiplying the right hand-side of the bubble problems with a constant. We reported the optimal values of this constant for  equilateral triangular elements. Various numerical experiments proved the robustness of the AB method in terms of the parameters provided. The AB method is able to solve the Helmholtz problem in 2D up to $ch=3.5$, efficiently. The numerical tests showed that the AB method is by far better than the pseudo-adaptive bubbles method and the fourth order method. We provided analysis to prove that the AB method mitigates the pollution error substantially.

\bibliographystyle{plain}
\bibliography{ref}






% \end{thebibliography}

\end{document}